\documentclass[10pt]{article}
\usepackage{amsmath,amsfonts,amssymb,amsthm,graphicx}
\usepackage{enumerate}
\usepackage{etoolbox}
\usepackage{color}
\usepackage[square]{natbib}
\usepackage{diagbox}
\usepackage{algorithm}
\usepackage[noend]{algpseudocode}
\algrenewcommand\algorithmicrequire{\textbf{Input:}}
\usepackage[colorlinks=true,citecolor=blue,pdfpagemode=UseNone,pdfstartview=FitH]{hyperref}

\allowdisplaybreaks[4]

\newcommand{\const}{\textrm{const}}

\newcommand{\p}{\textit{p}}
\newcommand{\e}{\textit{e}}

\newcommand{\pN}{\Box_{\p}}
\newcommand{\eN}{\Box_{\e}}
\newcommand{\pP}{\Diamond_{\p}}
\newcommand{\eP}{\Diamond_{\e}}

\renewcommand{\d}{\,\mathrm{d}}
\newcommand{\dd}{\mathrm{d}}
           
\newcommand{\E}{\mathbb{E}}
\renewcommand{\P}{\mathbb{P}}
\newcommand{\R}{\mathbb{R}}

\newcommand{\id}{1}

\newcommand{\AAA}{\mathcal{A}}

\newcommand{\UUU}{\mathcal{U}}
\newcommand{\XXX}{\mathcal{X}}

\DeclareMathOperator{\AV}{AV}
\DeclareMathOperator{\AM}{AM}
\DeclareMathOperator{\BM}{BM}
\DeclareMathOperator{\VS}{VS}
\renewcommand{\complement}{\textsf{c}}

\newcommand{\BRCAI}{\textrm{BRCA1}}

\theoremstyle{plain}
\newtheorem{theorem}{Theorem}[section]

\newtheorem{proposition}[theorem]{Proposition}

\theoremstyle{definition}

\theoremstyle{remark}
\newtheorem{remark}[theorem]{Remark}

\newcommand{\citeyr}[1]{[\citeyear{#1}]}

\begin{document}
\title{Confidence and discoveries with \e-values}
\author{Vladimir Vovk\thanks%
  {Department of Computer Science,
  Royal Holloway, University of London,
  Egham, Surrey, UK.
  E-mail: \href{mailto:v.vovk@rhul.ac.uk}{v.vovk@rhul.ac.uk}.}
\and Ruodu Wang\thanks%
  {Department of Statistics and Actuarial Science,
  University of Waterloo,
  Waterloo, Ontario, Canada.
  E-mail: \href{mailto:wang@uwaterloo.ca}{wang@uwaterloo.ca}.}}

\maketitle
\begin{abstract}
  \smallskip
  We discuss systematically two versions of confidence regions:
  those based on \p-values and those based on \e-values, a recent alternative to \p-values.
  Both versions can be applied to multiple hypothesis testing,
  and in this paper we are interested in procedures that control the number of false discoveries
  under arbitrary dependence between the base \p- or \e-values.
  We introduce a procedure that is based on \e-values
  and show that it is efficient both computationally and statistically
  using simulated and real-world datasets.
  Comparison with the corresponding standard procedure based on \p-values is not straightforward,
  but there are indications that the new one performs significantly better in some situations.

  \bigskip

  \noindent
  The journal version of this paper is to appear in \emph{Statistical Science}.
  This version has been further revised.
  For the most up-to-date version,
    see \href{http://alrw.net/e/}{http://alrw.net/e (Working Paper 3)}.
\end{abstract}

\section{Introduction}
\label{sec:introduction}

Starting from the introduction of confidence regions in the work of Jerzy Neyman \citeyr{Neyman:1934},
confidence estimation and hypothesis testing have been regarded as dual tasks.
We start our discussion from hypothesis testing and then extend it to confidence estimation.

The usual approaches to hypothesis testing and confidence estimation are based on \p-values,
but our emphasis will be on alternative approaches based on \e-values,
as discussed in, e.g., \citet{Shafer:2021} (who uses ``betting score'' for our ``\e-value''),
\citet[Section 11.5]{Shafer/Vovk:2019} (who use ``Skeptic's capital''),
\citet{Grunwald/etal:arXiv1906}, and \citet{Vovk/Wang:2021} (who proposed the term ``\e-values'').

\textit{E}-values can be defined as values taken by \e-variables,
and an \emph{\e-variable} is a random variable taking values in $[0,\infty]$
whose expectation is at most 1 under the null hypothesis.
In many areas of statistics \e-variables appear naturally as likelihood ratios:
if $Q$ is a simple null hypothesis and $Q'$ is an alternative probability measure,
the Radon--Nikodym derivative $\d Q'/\d Q$ is an \e-variable.
In Bayesian statistics, $Q$ or $Q'$ or both may be defined as marginal probability measures
for Bayesian models,
in which case likelihood ratios are known as Bayes factors.
The fundamental monograph treating Bayes factors is Jeffreys's \citeyr{Jeffreys:1961};
see, e.g., \citet{Ly/etal:2016} for a recent appreciation.
The notions of \e-values and Bayes factors coincide for simple null hypotheses but diverge for composite ones
(for \e-variables, the expectation should be at most 1 under any probability measure in the null hypothesis).

The existing statistical methods are often divided into Bayesian and classical
(we will say more about the latter in Section \ref{sec:comparison}).
While \p-values are the standard classical tool of hypothesis testing,
Bayes factors are the standard Bayesian tool \citep{Benjamini/etal:2021}.
One way of looking at \e-values is as a way of modelling Bayes factors
inside classical statistics inasmuch as they do not require prior distributions
for their definition.
This hints at the difficulty of comparisons between results
based on \p-values and those based on \e-values;
it is a manifestation of the oft-acknowledged chasm between classical and Bayesian statistics.

Roughly,
the Bayesian interpretation of an \e-variable $\d Q'/\d Q$ is that,
when deciding between $Q'$ and $Q$ as possible explanations for the data
and observing a very large \e-value,
the optimal decision is to reject $Q$
unless the prior probability of $Q$ is high
or a mistaken rejection of $Q$ is much more costly
than a mistaken rejection of $Q'$
(see \citet[Proposition 6.1]{Bernardo/Smith:2000}
for a precise decision-theoretic statement).

Another important and popular source of \e-values,
especially in the context of sequential observations,
is \emph{\e-processes}, which are stochastic processes $(E_t)_{t\ge 0}$
such that $E_\tau$ is an \e-variable for any stopping time $\tau$ (with respect to a pre-specified filtration);
see, e.g., \citet{Shafer/Vovk:2019}, \citet{Grunwald/etal:arXiv1906}, \citet{Vovk/Wang:2021}, and \citet{Wang/Ramdas:2022}.
The use of \e-processes ensures validity under optional stopping
and allows sequential update of statistical evidence.
These advantages are discussed extensively in the existing literature and are not the focus of this paper.

If $E$ is an \e-variable and $\alpha>1$,
Markov's inequality implies that, under the null hypothesis,
$E\ge\alpha$ with probability at most $1/\alpha$.
Therefore, observing a large value of $E$ provides evidence against the null hypothesis
in classical statistics as well.
In typical uses of \e-values, however, we are not given a threshold $\alpha$ in advance (or ever),
and simply regard an \e-value as the strength of evidence against the null.

The defining property of a \p-value is that it is $\alpha$ or less with a probability of at most $\alpha$.
This definition involves a quantifier over thresholds $\alpha$ and sometimes is considered misleading
in situations where no threshold $\alpha$ is fixed in advance.
There have been proposals to turn (``calibrate'') \p-values into Bayes factors \citep{Sellke/etal:2001} to help intuition,
and \citet[Appendix B]{Jeffreys:1961} proposes an informal correspondence between \p-values and Bayes factors;
both can be used for establishing connections between \p-values and \e-values.
Ways of turning \p-values into \e-values and vice versa
have been systematically discussed in \citet[Section~2]{Vovk/Wang:2021}
and are the topic of Section~\ref{sec:comparison}.
They provide ways of comparing results based on \e-values and \p-values,
albeit crude ones.

An area of statistics where we can see both \e-values and \p-values in action
is controlling the number of false discoveries in multiple hypothesis testing.
A known procedure of controlling the number of false discoveries
\citep{Genovese/Wasserman:2004,Goeman/Solari:2011local,Goeman/etal:2019Biometrika},
which we call the \emph{GWGS procedure},
uses \p-values, but it can be easily adapted to \e-values.
In this paper we demonstrate the performance of both versions.

Both versions of the GWGS procedure control the number of false discoveries in a stronger sense
than the well-known procedure of \citet{Benjamini/Hochberg:1995} controlling the false discovery rate (FDR).
Whereas FDR is the expected value of the false discovery proportion,
the GWGS procedure provides upper confidence bounds on the number of false discoveries.
Procedures that control FDR using \e-values are studied by \citet{Wang/Ramdas:2022}.

The GWGS procedure involves, at least implicitly,
merging several \p-values into a single \p-value.
Merging \p-values is difficult: see, e.g., \citet{Vovk/Wang:2020,Vovk/Wang/Wang:2022}.
The situation with \e-values is radically different:
arithmetic averaging is essentially the only symmetric method of merging
\citep[Proposition 3.1]{Vovk/Wang:2021}.
This contrast shows in the observation that the \e-version of the GWGS procedure
produces seemingly better results than the \p-version;
we cannot be more categorical since comparison between \p-values and \e-values
is not straightforward.

We start the main part of the paper by discussing testing in Sections~\ref{sec:testing}
and~\ref{sec:comparison},
defining confidence regions in Section \ref{sec:confidence_regions},
and repackaging them as necessity measures in Section~\ref{sec:confidence_measures}.
In Section~\ref{sec:control} we introduce the \e-version of the GWGS procedure,
postponing the \p-version to an appendix.
A special case that is easy to visualize is introduced
under the name of discovery \e-matrices.
In Sections~\ref{sec:simulation} and~\ref{sec:empirical}
we demonstrate the advantages of the \e-version
in simulation and empirical studies, respectively. 
In Section~\ref{app:AM} we give its computationally efficient implementation.
Section~\ref{sec:conclusion} concludes.

The main content of the paper is complemented by five appendixes, \ref{app:sigma}--\ref{app:truth}.
Appendix~\ref{app:sigma} contains some further information on a toy example in the main paper.
Appendix~\ref{app:general} explores other procedures
of controlling false discoveries with \e-values,
including the one based on a Bonferroni-type procedure of merging \e-values.
If the goal is family-wise validity,
such procedures (including the one in \citet{Holm:1979}) usually work very well,
but if the goal is to control the number of false discoveries,
they work much worse than arithmetic averaging.
In this appendix we also discuss a Simes-type procedure based on \e-values.
Appendix~\ref{app:GWGS} makes connections
with Goeman and Solari's \citeyr{Goeman/Solari:2011local} work explicit.
As \citet[Supplementary material]{Hemerik/etal:2019} explain,
the method of \citet{Goeman/Solari:2011local} is equivalent
to a method in \citet{Genovese/Wasserman:2004}.
Appendix~\ref{app:boosting} points out the importance of generalized Bayes factors.
Finally, Appendix~\ref{app:truth} summarizes results of further biomedical studies
related to the dataset that we use in Section~\ref{sec:empirical}.

\section{Three approaches to hypothesis testing}
\label{sec:testing}

\begin{figure}[htbp]
\begin{center}
  \includegraphics[width=0.4\textwidth]{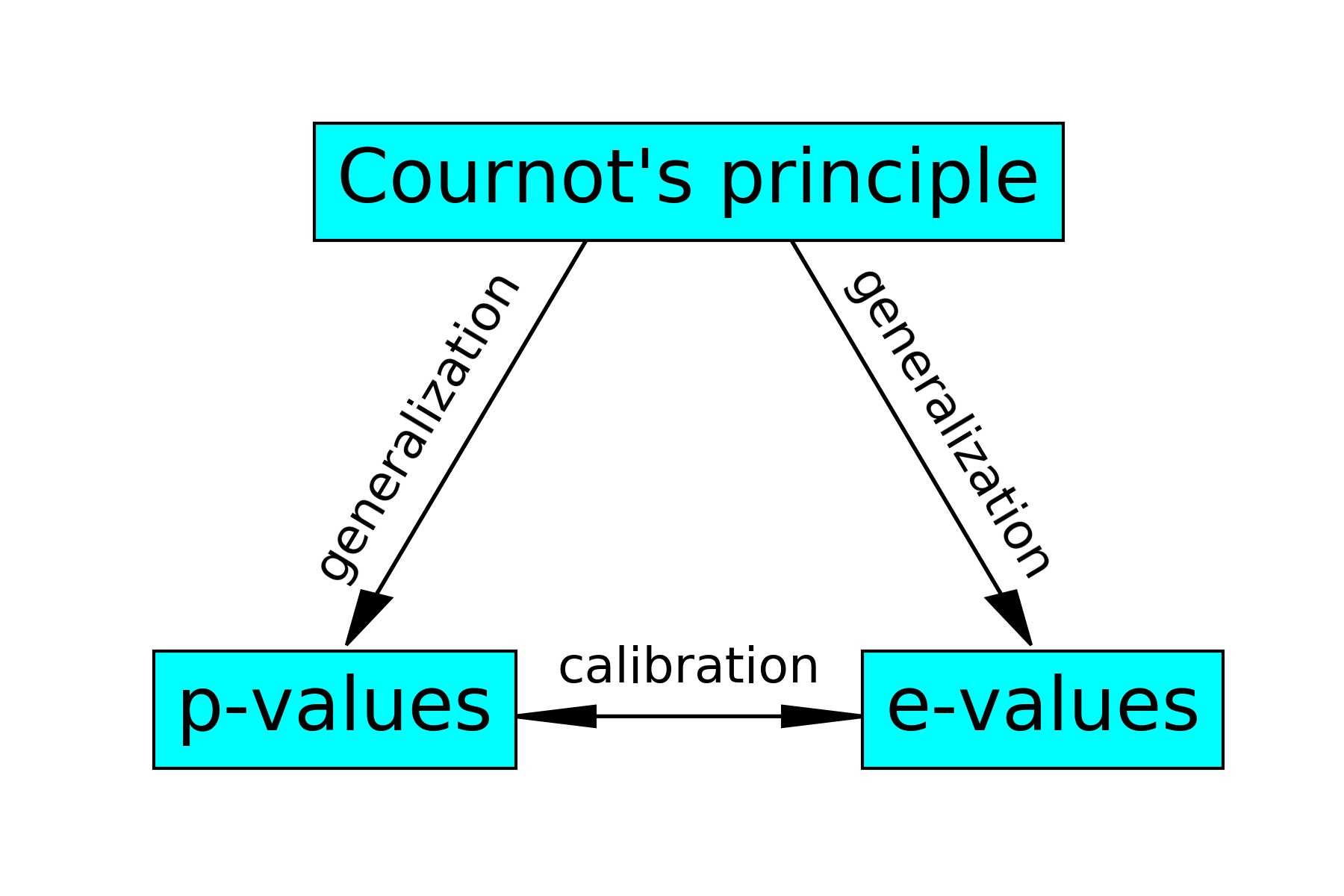}
\end{center} 
\caption{Cournot's principle and its two generalizations}\label{fig:Cournot}
\end{figure}

The basic principle of hypothesis testing is sometimes referred to
as Cournot's principle \citep{Shafer:2007-short}.
Augustin Cournot's bridge between probability theory and the world
is that if a given event has a small probability, we do not expect it to happen.
It is shown at the top of Figure~\ref{fig:Cournot}
and has entered (without its name) countless statistics textbooks:
the simplest approach to hypothesis testing consists in selecting \emph{a priori}
a critical region $A$ of a small probability under the null hypothesis
and rejecting the null hypothesis when $A$ happens.
Cournot's principle is the basis of the classical approach to statistics;
it was known to and used by James Bernoulli \citeyr{Bernoulli:1713},
and Cournot's \citeyr{Cournot:1843} contribution was to say
that this is the \emph{only} bridge.

We are mostly interested in two generalizations of Cournot's principle.
To give formal definitions,
we fix a measurable space $(\Omega,\AAA)$.
This is our \emph{sample space};
to complete it to a probability space we need a probability measure $Q\in\mathfrak{P}(\Omega)$,
where $\mathfrak{P}(\Omega)$ is the set of all probability measures on $(\Omega,\AAA)$.

A \emph{statistical model} is a family $(Q_{\theta}\mid\theta\in\Theta)$ of probability measures on $(\Omega,\AAA)$.
We do not require measurability, in any sense, of $Q_{\theta}$ in $\theta$;
in particular, the \emph{parameter space} $\Theta$ is just a set (not a measurable space).
We are mostly interested in the case where $\Theta=\mathfrak{P}(\Omega)$
and $Q_{\theta}=\theta$ for all $\theta\in\Theta$,
but in the first few sections our exposition will be general,
which may make our definitions more familiar to some of the readers.
For a given parameter $\theta\in\Theta$,
we have the notion of expectation $\E^\theta(E):=\int E \d Q_\theta \in [0,\infty]$
for each extended random variable $E$ taking nonnegative values
(we call it ``extended'' since it may take value $\infty$)
and the notion of probability $\P^\theta(A):=\E^\theta(\id_A)=Q_\theta(A)$ for each event $A\in\AAA$.

A \emph{simple statistical hypothesis} is an element $\theta$ of $\Theta$.
A \emph{statistical hypothesis} (or \emph{composite statistical hypothesis}, or simply \emph{hypothesis})
is a set $H\subseteq\Theta$ of parameters.
We embed the simple statistical hypotheses into the composite statistical hypotheses
by identifying $\theta\in\Theta$ with the corresponding singleton $\{\theta\}\subseteq\Theta$.
We will say ``null hypothesis'' to emphasize that we are interested
in whether the hypothesis should be rejected in view of the data $\omega\in\Omega$.

We already mentioned that the most basic way of testing a simple hypothesis $\theta\in\Theta$
is to choose a critical region $A\in\AAA$ with probability $\P^\theta(A)\le\alpha$,
$\alpha$ (the \emph{size}) being a small positive number,
and to reject the hypothesis $\theta$ at level $\alpha$ after observing an outcome $\omega\in A$.
A disadvantage of this way of testing is that it is binary;
either we completely reject the null hypothesis or we find no evidence whatsoever against it.
We will discuss two ways to graduate the notion of a critical region:
the classical one using \p-values and a more recent one using \e-values.

A \emph{\p-variable} for testing a simple hypothesis $\theta$
is a nonnegative random variable $P$ such that, for any $\alpha\in(0,1)$,
$\P^\theta(P\le\alpha)\le\alpha$.
For each threshold $\alpha$ we have a critical region $\{P\le\alpha\}$,
and a \p-variable provides a nested family of critical regions.
An \emph{\e-variable} for testing a simple hypothesis $\theta\in\Theta$
is a nonnegative extended random variable $E$ such that $\E^\theta(E)\le1$.

Suppose we are testing a simple null hypothesis $\theta$
(it might correspond to a default parameter value).
In \p-testing, we choose a \p-variable $P$ in advance
and reject the null hypothesis $\theta$
when the observed value $P(\omega)$ of $P$ (the \emph{p-value}) is small,
and in \e-testing, we choose an \e-variable $E$ in advance
and reject the null hypothesis $\theta$
when the observed value $E(\omega)$ of $E$ (the \emph{e-value}) is large.
In both cases,
we get a measure of the amount of evidence found against the null hypothesis.

We can embed basic testing into both \p-testing and \e-testing:
namely, to each critical region $A$ corresponds the \p-variable
\begin{equation}\label{eq:basic_to_p}
  P(\omega)
  :=
  \begin{cases}
    \alpha & \text{if $\omega\in A$}\\
    1 & \text{if not}
  \end{cases}
\end{equation}
and \e-variable
\begin{equation}\label{eq:basic_to_e}
  E(\omega)
  :=
  \begin{cases}
    1/\alpha & \text{if $\omega\in A$}\\
    0 & \text{if not},
  \end{cases}
\end{equation}
where $\alpha$ is the size of the critical region $A$.
These two random variables carry the same information as $A$.
This justifies the two arrows marked ``generalization'' in Figure~\ref{fig:Cournot}.

The special case of basic testing corresponds
to concentrating on only one threshold,
denoted $\alpha$ in the case of \p-values, \eqref{eq:basic_to_p},
and $1/\alpha$ in the case of \e-values, \eqref{eq:basic_to_e}.
It is instructive to see how we could extend
the basic \p-variable \eqref{eq:basic_to_p}
and the basic \e-variable \eqref{eq:basic_to_e}.
It is easy to extend \eqref{eq:basic_to_p};
e.g., we can take another critical region $A'\supset A$
of size $\alpha'>\alpha$ and define a p-variable,
\begin{equation*}
  P(\omega)
  :=
  \begin{cases}
    \alpha & \text{if $\omega\in A$}\\
    \alpha' & \text{if $\omega\in A'\setminus A$}\\
    1 & \text{if $\omega\in\Omega\setminus A'$},
  \end{cases}
\end{equation*}
that strongly dominates \eqref{eq:basic_to_p}.
As it were, for each $\alpha$ we have a separate budget of $\alpha$
that can be spent on a critical region.
On the other hand, there is no way to improve the \e-variable \eqref{eq:basic_to_e}
in a non-trivial way (make it larger on a set of positive probability).
Now we have a single budget of 1,
which has been fully spent in \eqref{eq:basic_to_e}.

The definitions of critical regions, \p-variables, and \e-variables
extend to the case of composite hypotheses as follows.
A \emph{critical region} of size $\alpha$ for a composite hypothesis $H$
is an event $A\in\AAA$ satisfying $\P^\theta(A)\le\alpha$ for all $\theta\in H$.
A \emph{\p-variable} for testing a composite hypothesis $H$
is a nonnegative random variable $P$ such that, for any $\alpha\in(0,1)$,
$\P^\theta(P\le\alpha)\le\alpha$ for all $\theta\in H$.
And an \emph{\e-variable} for testing a composite hypothesis $H$
is an extended nonnegative $E$
satisfying $\E^\theta(E)\le1$ for all $\theta\in H$.

While \p-variables (referred to as valid \p-values in \citet[Definition 8.3.26]{Casella/Berger:2002})
are standard, \e-variables \citep{Shafer:2021,Vovk/Wang:2021} have not been used widely.

Observing a small \p-value or a large \e-value provide evidence against $H$.
It is convenient to have conventional thresholds for \p-values and \e-values.
For \p-values, the standard thresholds are 1\% and 5\%,
and they go back to Fisher.
If $p\le0.05$, we say that the evidence against the null hypothesis is significant,
and if $p\le0.01$, we say that the evidence is highly significant.
For \e-values, we will use Jeffreys's [\citeyear[Appendix B]{Jeffreys:1961}] rule of thumb:
\begin{itemize}
\item
  If the \e-value is below 1,
  the null hypothesis is supported.
  In our plots (such as in Figure~\ref{fig:small} below)
  in the experimental sections, \ref{sec:simulation} and \ref{sec:empirical},
  such \e-values will be shown in dark green.
\item
  If the \e-value is in the interval $(1,\sqrt{10})\approx(1,3.16)$,
  the evidence against the null hypothesis is not worth more than a bare mention.
  Such \e-values will be shown in green.
\item
  If the \e-value is in $(\sqrt{10},10)\approx(3.16,10)$,
  the evidence against the null hypothesis is substantial.
  Shown in yellow.
\item
  If it is in $(10,10^{3/2})\approx(10,31.6)$,
  the evidence against the null hypothesis is strong.
  Shown in red.
\item
  If it is in $(10^{3/2},100)\approx(31.6,100)$,
  the evidence against the null hypothesis is very strong.
  Shown in dark red.
\item
  If the \e-value exceeds $100$,
  the evidence is decisive.
  Shown in black.
\end{itemize}

\subsection*{Fisher's and Neyman--Pearson's views of testing}

A common view is that, in our terminology,
Fisher preferred \p-testing, whereas Neyman preferred basic testing
(with the null hypothesis complemented by an alternative hypothesis).
The full story is, however, more complex: see, e.g., \citet[Section 4.4]{Lehmann:2011}.

Fisher's interpretation of hypothesis testing was in terms of a disjunction
\citep[Section III.1]{Fisher:1973-short}.
If $A$ is a critical region of a small size $\alpha$ and we observe an outcome in $A$,
then \emph{either the null hypothesis is wrong or ``a rare chance has occurred''}.
To avoid any frequentist connotations,
we may express it in the equivalent form
\emph{the null hypothesis is wrong unless the outcome is strange}
(``unless'' being one of the ways to express the idea of disjunction
\citep[Section 14]{Kleene:1967}).
A similar interpretation is applicable to \p-values and \e-values:
e.g., if we observe a large \e-value,
then the null hypothesis is wrong unless the outcome is strange.

In Neyman and Pearson's approach to hypothesis testing,
a big role is played by alternative hypotheses.
In \e-testing, the notion of an alternative hypothesis plays a less independent role:
choosing an \e-variable $E$ can often be interpreted as choosing an alternative hypothesis
in such a way that $E$ is the likelihood ratio of the alternative hypothesis to the null
\citep[2.2]{Shafer:2021}.

\subsection*{A toy example}

\begin{figure}
  \begin{center}
    \includegraphics[width=0.48\textwidth]{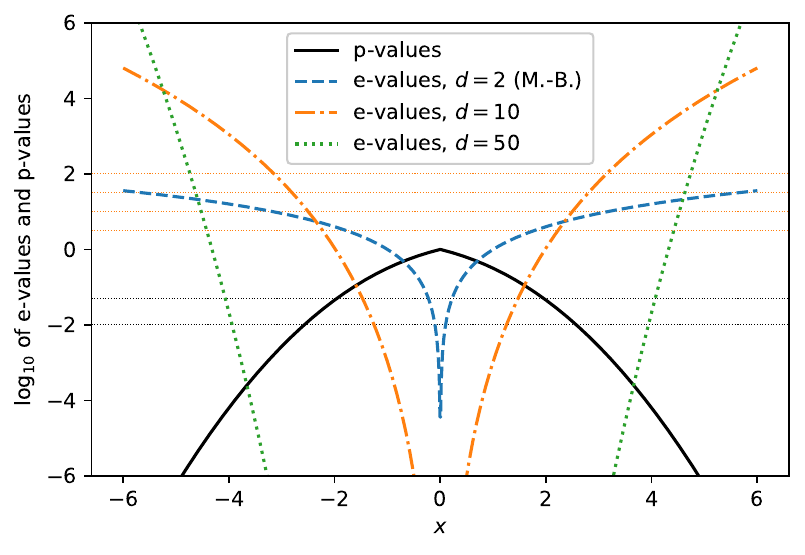}
    \includegraphics[width=0.48\textwidth]{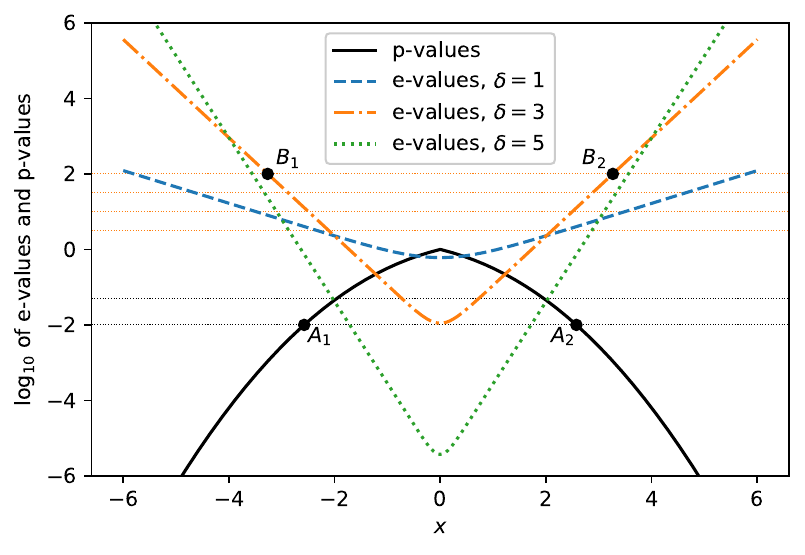}
  \end{center}
  \caption{The \p-values (black solid lines) and \e-values on the decimal log scale
    for testing the null hypothesis $N(0,1)$.
    Left panel: for the signed $\chi$ alternatives.
    Right panel: for the alternatives in the family $N(\mu,1)$.}
  \label{fig:toy}
\end{figure}

Let us see how these definitions work in a simple example.
We would like to test the null hypothesis $x\sim N(0,1)$
given an observation $x\in\Omega:=\R$.
Suppose we believe that $\lvert x\rvert$ reflects
the amount of evidence against the null hypothesis.
Therefore, we will be interested in \p-variables $P(x)$ and \e-variables $E(x)$
that depend on $x$ only via $\lvert x\rvert$
and are monotonic functions
(increasing for $E$ and decreasing for $P$) of $\lvert x\rvert$.
There is a unique  \p-variable $P$ (uniformly distributed on $[0,1]$) satisfying this property,
namely $P(x):=2\Phi(-\lvert x\rvert)$,
where $\Phi$ is the standard Gaussian distribution function.
On the other hand, there is a huge variety of \e-variables
satisfying this property.
A natural class of such \e-variables is
\begin{equation}\label{eq:d-family}
  E(x)
  :=
  \frac{\lvert x\rvert^d}{\pi^{-1/2}2^{d/2}\Gamma\left(\frac{d+1}{2}\right)},
  \quad
  d>0,
\end{equation}
where the denominator is just the normalizing constant
ensuring $\int E\d N(0,1)=1$
(the $d$th absolute moment of the standard Gaussian distribution,
which is well known and easily found by direct integration).

Figure~\ref{fig:toy} gives the \p-values as the black solid line
in both panels
and gives the \e-values for $d\in\{2,10,50\}$ in the left panel.
On Fisher's scale, \p-values are significant
when their decimal logarithms drop below $\log_{10}0.05$
and highly significant when they drop below $-2$;
these levels are shown as thin black lines.
Jeffreys's levels $0.5$, $1$, $1.5$, and $2$ for \e-values on the $\log_{10}$ scale
are shown as thin orange lines.
The \e-variables are not comparable, in the sense that none of them
dominates any other everywhere.

Each of the \e-variables in the left panel of Figure~\ref{fig:toy}
defines an alternative to the null hypothesis, as discussed above,
so that the \e-variable becomes
the likelihood ratio of the alternative to the null hypothesis.
The alternative hypothesis corresponding to \eqref{eq:d-family}
has the density proportional to $\lvert x\rvert^d\exp(-x^2/2)$.
Since $x$ ranges over $\R$, it is not exactly the $\chi$ density
with $d+1$ degrees of freedom, which we will denote $\chi_{d+1}$,
but it is a slight variation:
after generating $x$ from $\chi_{d+1}$, we change its sign
(i.e., multiply it by $-1$)
with probability $1/2$.
We will call this alternative distribution the \emph{signed $\chi_{d+1}$ distribution};
for $d=2$ this is the \emph{signed Maxwell--Boltzmann distribution},
which is abbreviated to ``M.-B.''
in the legend in the left panel of Figure~\ref{fig:toy}.

The signed $\chi$ alternatives,
corresponding to $E(x)\propto\lvert x\rvert^d$,
appear to be  the simplest unconstrained choice,
but a more standard approach is to look for alternatives
inside a parametric family of distributions.
Let us embed $N(0,1)$ into the statistical model $N(\mu,1)$, $\mu\in\R$
(there are other natural embeddings,
and in Appendix~\ref{app:sigma} we will also discuss
the embedding into the statistical model $N(0,\sigma^2)$, $\sigma>0$).
The right panel of Figure~\ref{fig:toy} shows three more \e-variables,
which are based on the likelihood ratios
\begin{equation}\label{eq:E_mu}
  E^{(\delta)}(x)
  :=
  \frac{\dd N(\delta,1)}{\dd N(0,1)}(x)
  =
  \frac{\exp(-(x-\delta)^2/2)}{\exp(-x^2/2)}
  =
  \exp(\delta x - \delta^2/2).
\end{equation}
To obtain an \e-variable that depends on $x$ only via $\lvert x\rvert$,
let us average $E^{(\delta)}$ and $E^{(-\delta)}$:
\begin{equation}\label{eq:average}
  \bar E^{(\delta)}(x)
  :=
  \frac{E^{(\delta)}(x)+E^{(-\delta)}(x)}{2}.
\end{equation}
The right panel of Figure~\ref{fig:toy} shows
$\bar E^{(1)}$, $\bar E^{(3)}$, and $\bar E^{(5)}$.

Four points are indicated in the right panel
to illustrate the interpretation of our plots.
The points $A_1$ and $A_2$ are at the intersection of the graph of the \p-variable
with the horizontal line at level $-2$,
and the points $B_1$ and $B_2$ are
at the intersection of the graph of the \e-variable with $\delta:=3$
with the horizontal line at level $2$.
The $x$-coordinates of the points $A_1$ and $A_2$ are approximately $\pm2.58$,
and the $x$-coordinates of the points $B_1$ and $B_2$ are approximately $\pm3.27$.
The observations $x$ with $\lvert x\rvert$ exceeding (approximately) $3.27$
provide decisive evidence against the null hypothesis $N(0,1)$,
and the other observations do not,
according to Jeffreys's scale.
Similarly, the observations $x$ with $\lvert x\rvert$ exceeding (approximately) $2.58$
provide highly significant evidence against the null hypothesis $N(0,1)$,
while the other observations do not,
according to Fisher's scale.
Using similar interpretation for the left panel,
we can see, e.g., that the \e-variable for $d=2$ makes a wider range of observations $x$
provide substantial evidence against the null hypothesis
than the \e-variable for $d=10$ does,
whereas for strong evidence we have the opposite situation.

Later in the paper we will use methods related
to both panels of Figure~\ref{fig:toy}.
In our simulation studies in Section~\ref{sec:simulation},
we will generate observations from $N(\mu,1)$
and use again the likelihood ratios \eqref{eq:E_mu}
that we used in the right panel.
In our empirical studies in Section~\ref{sec:empirical},
where we have no idea of the true distribution of the data,
we use the \e-variables proportional to $\lvert x\rvert^d$,
as in the left panel.

\section{Is it possible to compare \e-values and \p-values?}
\label{sec:comparison}

Starting from the next section we will describe various methods
based on \p-values and \e-values.
In Section~\ref{sec:introduction} we already alluded to difficulties
of comparing such results.
There is no overarching testing framework (at least at this time)
containing both \p-testing and \e-testing
that could be used for comparing such results.
The best we can do rigorously is to convert, albeit imperfectly,
\p-values to \e-values and vice versa.

Sometimes the user of statistical procedures
has a clear preference for \p-values or \e-values.
These are some possible categories of users
(this is not an exhaustive list, of course):
\begin{enumerate}
\item
  Some users will find the frequentist interpretation of \p-variables $P$ appealing:
  for any threshold $\alpha$, the long-run frequency of observing $P\le\alpha$
  in a sequence of independent identical trials is at most $\alpha$.
  This is typical of the frequentist school of classical statistics.
\item
  Other users will prefer a direct application of Cournot's principle:
  for a small $\alpha$ and pre-specified $P$, we do not expect to observe $P\le\alpha$
  under the null hypothesis.
  This school is referred to as Bernoullian statistics
  by Glenn Shafer \citep{Shafer:BFF},
  following Francis Edgeworth, Richard von Mises, Arthur Dempster, and Ian Hacking.
\item
  Another category, representing Bayesian statistics,
  will like the Bayesian interpretation of \e-variables
  referred to in Section~\ref{sec:introduction}.
\item
  Our final category will accept the betting interpretation of an \e-variable $E$
  (see, e.g., \citet{Shafer:2021}):
  the \e-value $E(\omega)$ is the pay-off of a lottery
  that is fair under the null hypothesis,
  and a large pay-off casts doubt on the null hypothesis.
  The idea of betting is often regarded as an important ingredient of Bayesian statistics
  (see, e.g, \citet[3.3.5]{deFinetti:2017}),
  but it is used there in a very different way,
  in the form of no Dutch book requirement.
\end{enumerate}
Communication may be easier between the users in the first two categories (classical statistics),
or between the users in the last two categories.
But otherwise, we need methods of conversion between \p- and \e-values.
Therefore, as a first step we discuss rigorous ways of turning \p-values into \e-values
(known as calibrating \p-values)
and vice versa.
For further details, see \citet[Section~2]{Vovk/Wang:2021}.

A decreasing function $f:[0,1]\to[0,\infty]$ is a \emph{calibrator}
if, for any \p-variable $P$, $f(P)$ is an \e-variable.
In other words, a calibrator transforms \p-values to \e-values.
A very natural family of calibrators is
\begin{equation}\label{eq:calibrator}
  f_{\kappa}(p)
  :=
  \kappa p^{\kappa-1},
\end{equation}
where $\kappa\in(0,1]$.
The maximum possible \e-value
\begin{equation}\label{eq:VS}
  \VS(p)
  :=
  \max_{\kappa\in(0,1]}
  f_{\kappa}(p)
  =
  \begin{cases}
    -\exp(-1)/(p\ln p) & \text{if $p\le\exp(-1)$}\\
    1 & \text{otherwise}
  \end{cases},
  \quad
  p\in(0,1],
\end{equation}
attainable by this family will be referred to as the \emph{VS bound}
(abbreviating ``Vovk--Sellke bound'' \citep[Section 11.5]{Sellke/etal:2001,Shafer/Vovk:2019}),
but due to the maximum operation, $\VS(P)$ need not be an \e-variable
even if $P$ is a \p-variable.

In the opposite direction,
a decreasing function $f:[0,\infty]\to[0,1]$ is an \emph{\e-to-\p\ calibrator}
if, for any \e-variable $E$, $f(E)$ is a \p-variable.
It is a function transforming \e-values to \p-values.
As explained and formalized in \citet[Proposition~2.2]{Vovk/Wang:2021},
\begin{equation}\label{eq:e-to-p_calibrator}
  t\in[0,\infty]\mapsto\min(1,1/t)
\end{equation}
is the only reasonable \e-to-\p\ calibrator.

In general, calibrating \p-values and \e-values are crude processes.
A strong \e-value of $20$ barely attains statistical significance
when transformed into a \p-value (namely, $5\%$) using \eqref{eq:e-to-p_calibrator}.
The VS bound for the borderline significant \p-value of $5\%$ is approximately $2.456$,
and so ``is not worth more than a bare mention'', according to Jeffreys.
The low ``round-trip efficiency'' in the domain of \p-values can be illustrated by
\begin{equation}\label{eq:round_trip}
  1/\VS(0.005) \approx 0.072.
\end{equation}
The round trip turns the highly significant \p-value of $0.5\%$ into the non-significant \p-value of $7.2\%$.
And this is despite the VS bound being achievable as \e-value only in hindsight.

In view of the low round-trip efficiency,
it is natural to expect that users of statistical procedures
who insist on using \p-values will be best served by methods producing directly \p-values.
A method producing \e-values will have to be vastly superior
to result in better, or even equally good, \p-values after conversion.
Interestingly, we will see such an example in Section~\ref{sec:simulation}
(see the discussion of Figure~\ref{fig:big}).
Symmetrically, a method producing \p-values will have to be vastly superior
to a method producing directly \e-values
in order to result in better or equally good \e-values after conversion.
The caveat here is that the result of comparison
still depends on using the bound~\eqref{eq:VS}
(using the \e-to-\p\ calibrator \eqref{eq:e-to-p_calibrator} is uncontroversial).

The user who is uncertain whether to use \p-values or \e-values
usually needs a more accurate comparison than that provided
by the crude procedures of calibration and \e-to-\p\ calibration.
We do not have objective ways of doing that.
One subjective way to compare results using \e-values to those using \p-values
is to appeal to Jeffreys's [\citeyear[Appendix B]{Jeffreys:1961}] authority:
``Users of these tests\label{p:Jeffreys} speak of the 5 per cent.\ point in much the same way
as I should speak of the $K = 10^{-1/2}$ point,
and of the 1 per cent.\ point as I should speak of the $K = 10^{-1}$ point.''
In our terminology, people doing \p-testing speak of a \p-value of $5\%$ (resp.\ $1\%$)
in much the same way as Jeffreys should speak of an \e-value of $10^{1/2}$ (resp.\ $10$).
The approximate equivalences are
\begin{equation}\label{eq:Jeffreys}
  \text{\p-value of $5\%$} \sim \text{\e-value of $10^{1/2}\approx3.16$}
  \qquad
  \text{\p-value of $1\%$} \sim \text{\e-value of $10$}.
\end{equation}

Another subjective way is to use Good's \citep[Appendix~IV]{Good:1958} rule of thumb.
According to Good,
the \e-value corresponding to a \p-value of $p$ should lie in the range
\begin{equation}\label{eq:Good_int}
  \left(
    \frac{1}{30p},
    \frac{3}{10p}
  \right)
\end{equation}
when $0.001<p<0.2$
(which Good felt were the values of $p$ that are usually of most practical interest).
In Good's picture the \p-value of $p$ is obtained using the standard recipe
from the Bayes factor as test statistic (this condition is always satisfied in this paper).
If we take the geometric mean $1/(10p)$ of the end-points of the interval \eqref{eq:Good_int},
we will obtain
\begin{equation}\label{eq:Good}
  \text{\p-value of $5\%$} \sim \text{\e-value of $2$}
  \qquad
  \text{\p-value of $1\%$} \sim \text{\e-value of $10$}
\end{equation}
in place of \eqref{eq:Jeffreys}.
While \eqref{eq:Jeffreys} and \eqref{eq:Good} are close,
Good acknowledges the significant uncertainty surrounding the correspondence.

A slightly more objective way of comparing methods based on \p-values and \e-values
is to consider their mathematical simplicity.
A great advantage of \e-values is that they are very easy to combine;
as we mentioned in Section~\ref{sec:introduction},
arithmetic averaging is essentially the only symmetric method of combination
\citep[Proposition 3.1]{Vovk/Wang:2021}.
This leads to simple and intuitive algorithms
(and in \citet{Vovk/Wang/Wang:2022} merging \e-values is even used
as a technical tool for designing admissible ways of merging \p-values).

\section{Confidence regions}
\label{sec:confidence_regions}

The notion of a confidence region was introduced by Neyman \citeyr{Neyman:1934,Neyman:1937}
only in its basic version.
(See \citet[Section 6.4]{Lehmann:2011} for Neyman's predecessors;
the word ``confidence'' is a translation of the Polish ``ufno\'s\'c'' \citep{Neyman:1941},
and Neyman's adjectival use of it was at first made fun of by his English listeners
\citep[comments by Bowley and Fisher]{Neyman:1934}.)
The \p-version is usually implicit, and the \e-version has not been used in mainstream statistics.
However, the \e-version has been used for a long time, in some form,
in the algorithmic theory of randomness \citep{Levin:1976uniform-Latin,Gacs:2005,Vovk/Vyugin:1993},
and in this paper we will use the terminology close to that of the algorithmic theory of randomness.

Let us fix a statistical model $(Q_{\theta}\mid\theta\in\Theta)$.
A \emph{basic test} of size $\alpha$ is a family of critical regions $(A_\theta\mid\theta\in\Theta)$ of size $\alpha$.
Therefore, for each simple statistical hypothesis $\theta$,
a basic test fixes a critical region $A_\theta$ for testing $\theta$:
$\P^\theta(A_\theta)\le\alpha$.

The interpretation of a basic test that is symmetric between the parameter space $\Theta$ and sample space $\Omega$
is that $\omega\in A_{\theta}$ means poor agreement between $\theta$ and $\omega$.
This binary relation of poor agreement and its complementary relation of good agreement have two sides:
\begin{itemize}
\item
  on the testing side,
  we start from $\theta$ and divide the $\omega$s into those that conform to $\theta$ ($\omega\notin A_\theta$)
  and those that do not ($\omega\in A_\theta$);
\item
  on the estimation side,
  we start from $\omega$ and divide the $\theta$s into those that agree with $\omega$ ($\omega\notin A_\theta$)
  and those that do not ($\omega\in A_\theta$).
\end{itemize}
In particular, on the estimation side we have the notion of a confidence estimator as introduced by Neyman
(cf.\ \citet[Figure IV]{Neyman:1934}):
the confidence estimator corresponding to a basic test $(A_\theta\mid\theta\in\Theta)$ is
\begin{equation}\label{eq:basic_estimator}
  \Gamma(\omega)
  :=
  \{\theta\in\Theta\mid\omega\notin A_\theta\}.
\end{equation}

In the context of a basic test $(A_\theta\mid\theta\in\Theta)$ of a small size $\alpha$
we may say that an outcome $\omega\in\Omega$ is \emph{strange} for a parameter value $\theta\in\Theta$
if $\omega\in A_\theta$.
According to Cournot's principle, we do not expect the outcome $\omega$ to be strange for the true $\theta$.
Our interpretation of the confidence region \eqref{eq:basic_estimator}
is that $\Gamma(\omega)$ covers the true $\theta$ unless $\omega$ is strange.

Graduated notions of a confidence estimator are discussed surprisingly rarely in statistics textbooks,
especially in full generality
(e.g., the popular textbook \citet{Cox/Hinkley:1974}
is one of the few places where they are discussed,
but only in the context of a linearly ordered parameter space $\Theta$).
A \emph{\p-test} is a family of \p-variables $(P_\theta\mid\theta\in\Theta)$,
and the corresponding \emph{\p-confidence regions} are defined as 
\begin{equation}\label{eq:p-estimator}
  \Gamma_{\p,\alpha}(\omega)
  :=
  \{\theta\in\Theta\mid P_\theta(\omega)>\alpha\},
  \quad
  \alpha\in(0,1).
\end{equation}
We regard $P_\theta(\omega)$ as a measure of agreement between $\theta$ and $\omega$,
with small values indicating poor agreement,
and define $\Gamma_{\p,\alpha}(\omega)$ to be the set of $\theta$
that agree with $\omega$ at level $\alpha$.
The definition of a p-confidence estimator
is only a slight variation on the definition of a basic estimator:
namely, \eqref{eq:p-estimator} can be obtained from \eqref{eq:basic_estimator}
by setting $A_\theta:=\{\omega\mid P_\theta(\omega)\le\alpha\}$
for each $\alpha$.
Notice that the  {\p-confidence regions} $\Gamma_{\p,\alpha}(\omega)$
are \emph{nested}:
$\alpha_1<\alpha_2$ implies
$\Gamma_{\p,\alpha_2}(\omega)\subseteq\Gamma_{\p,\alpha_1}(\omega)$;
this property is sometimes discussed or at least mentioned in statistics textbooks
(e.g., in \citet[Section 7.2]{Cox/Hinkley:1974},
\citet[9.5.2]{Casella/Berger:2002},
and \citet[Section 19.15]{Kendall:1999}).

Similarly, an \emph{\e-test} is a family of \e-variables $(E_\theta\mid\theta\in\Theta)$.
We also regard $E_\theta(\omega)$ as a measure of agreement between $\theta$ and $\omega$,
but now large values indicate poor agreement.
Analogously to \eqref{eq:p-estimator}, we define the \emph{\e-confidence regions} as
\begin{equation}\label{eq:e-estimator}
  \Gamma_{\e,\alpha}(\omega)
  :=
  \{\theta\in\Theta\mid E_\theta(\omega)<\alpha\},
  \quad
  \alpha\in(0,\infty).
\end{equation}
The definitions \eqref{eq:p-estimator} and \eqref{eq:e-estimator}
of \p-confidence regions and \e-confidence regions
generalize the basic definition \eqref{eq:basic_estimator},
which corresponds to using the \p-test and the \e-test
defined by \eqref{eq:basic_to_p} and \eqref{eq:basic_to_e},
respectively,
with added subscripts $\theta$.

The notions of \p-test and \e-test provide graduated notions of strangeness.
Let $\alpha>0$;
we will sometimes refer to it as the \emph{significance level}
(the interesting values are $\alpha<1$ for \p-testing and $\alpha>1$ for \e-testing).
In the context of a \p-test $(P_\theta\mid\theta\in\Theta)$,
we say that $\omega\in\Omega$ is \emph{$\alpha$-strange} for $\theta\in\Theta$
if $P_\theta(\omega)\le\alpha$
(i.e., if we reject $\theta$ at level $\alpha$ after observing $\omega$).
And in the context of an \e-test $(E_\theta\mid\theta\in\Theta)$,
we say that $\omega\in\Omega$ is \emph{$\alpha$-strange} for $\theta\in\Theta$
if $E_\theta(\omega)\ge\alpha$.
If there is any risk of confusion,
we will use the fuller expressions ``$(\text{\p},\alpha)$-strange'' and ``$(\text{\e},\alpha)$-strange''.

The interpretation of the confidence region \eqref{eq:p-estimator}
in terms of a Fisher-type disjunction
is that $\Gamma_{\p,\alpha}(\omega)$ covers the true $\theta$ unless $\omega$ is $(\p,\alpha)$-strange.
Similarly, we interpret \eqref{eq:e-estimator}
by saying that $\Gamma_{\e,\alpha}(\omega)$ covers the true $\theta$ unless $\omega$ is $(\e,\alpha)$-strange.

Starting from Section~\ref{sec:simulation},
we will visualize \e-confidence regions for a range of thresholds, including 10.
Inspired by the terminology of \citet[Appendix B]{Jeffreys:1961},
already discussed in Section~\ref{sec:testing},
we will refer to an \e-confidence region at level $1$ as \emph{weak},
at level $10^{1/2}$ as \emph{substantial},
at level $10$ as \emph{strong},
at level $10^{3/2}$ as \emph{very strong},
and at level $100$ as \emph{extremely strong}.

\subsection*{Simultaneous confidence regions}

Sometimes we are interested not in $\theta$ but in some derivative parameter
(as in Schervish's textbook [\citeyear[5.2.1]{Schervish:1995}]).
For example, if $\Omega=\R$ and $\Theta=\mathfrak{P}(\Omega)$,
we might be interested in the median of $\theta\in\Theta$.
Let $g:\Theta\to\Theta_g$ be the function mapping the original parameter $\theta$
to a new parameter, $g(\theta)$.

The confidence regions for the derived parameter $g(\theta)$ become:
\begin{equation}\label{eq:basic_estimator_g}
  \Gamma^g(\omega)
  :=
  \{g(\theta)\mid\theta\in\Theta\And\omega\notin A_\theta\}
\end{equation}
in place of \eqref{eq:basic_estimator},
\begin{equation*}
  \Gamma_{\p,\alpha}^g(\omega)
  :=
  \{g(\theta) \mid \theta\in\Theta \And P_\theta(\omega)>\alpha\},
  \quad
  \alpha\in(0,1),
\end{equation*}
in place of \eqref{eq:p-estimator}, and
\begin{equation}\label{eq:e-estimator_g}
  \Gamma_{\e,\alpha}^g(\omega)
  :=
  \{g(\theta) \mid \theta\in\Theta \And E_\theta(\omega)<\alpha\},
  \quad
  \alpha\in(0,\infty),
\end{equation}
in place of \eqref{eq:e-estimator}.

It is important that we can have a family of functions $g$,
and the confidence estimator \eqref{eq:basic_estimator_g} will be valid simultaneously for all of them,
provided the same basic test $(A_\theta\mid\theta\in\Theta)$ is used for all $g$.
The same is true for \p-confidence estimators and \e-confidence estimators;
what is important is that the notion of strangeness should not depend on $g$.
For example, for any family of functions $g:\Theta\to\Theta_g$,
the confidence region $\Gamma_{\e,\alpha}^g(\omega)$ in \eqref{eq:e-estimator_g}
contains $g(\theta)$ for all $g$ simultaneously unless the outcome is $\alpha$-strange
for the true parameter $\theta$.

\subsection*{Confidence regions in the toy example}

Here we continue our discussion of the toy example started in the previous section.
Now our statistical model $(Q_\theta\mid\theta\in\Theta)$ is $\Theta:=\R$
and $Q_\theta:=N(\theta,1)$ for all $\theta$.
For a fixed $\delta$, such as $\delta:=3$,
let us generalize \eqref{eq:average} to
\begin{equation}\label{eq:E-bar}
  \bar E^{(\delta)}_{\theta}(x)
  :=
  \frac{E^{(\delta)}_{\theta}(x)+E^{(-\delta)}_{\theta}(x)}{2},
\end{equation}
where, generalizing \eqref{eq:E_mu},
\begin{equation}\label{eq:E-one-sided}
  E^{(\delta)}_{\theta}(x)
  :=
  \frac{\dd N(\theta+\delta,1)}{\dd N(\theta,1)}(x)
  =
  \frac{\exp(-(x-\theta-\delta)^2/2)}{\exp(-(x-\theta)^2/2)}
  =
  \exp(\delta(x-\theta) - \delta^2/2).
\end{equation}
This gives us an \e-test.

Remember that the $x$-coordinate of the point $B_2$
in the right panel of Figure~\ref{fig:toy} is approximately $3.27$,
and let us fix $\delta:=3$.
Therefore, the extremely strong \e-confidence regions
(\e-confidence intervals in this case) are
\begin{equation}\label{e-confidence_1}
  \Gamma_{e,100}(x)
  \approx
  [x-3.27,x+3.27],
  \quad
  x\in\R,
\end{equation}
where ``$\approx$'' refers to $3.27$ being an approximate value.
For Jeffreys's other thresholds the \e-confidence intervals are
\begin{multline}\label{e-confidence_2}
  \Gamma_{e,10^{3/2}}(x)
  \approx
  [x-2.88,x+2.88],
  \quad
  \Gamma_{e,10}(x)
  \approx
  [x-2.50,x+2.50],\\
  \Gamma_{e,10^{1/2}}(x)
  \approx
  [x-2.11,x+2.11],
\end{multline}
and the \p-confidence intervals are, as usual,
$
  \Gamma_{p,0.01}(x)
  \approx
  [x-2.58,x+2.58]
$.

The \e-confidence intervals \eqref{e-confidence_1}--\eqref{e-confidence_2}
will change if the alternative hypotheses $\theta\pm3$ are replaced by other ones,
such as $\theta\pm1$ or $\theta\pm5$.
It can be considered an advantage of \p-confidence intervals,
and \p-values in general,
that for an important (albeit small) set of popular statistical models
there is no dependence on the choice of the alternative hypothesis.
This is closely related to the existence of uniformly most powerful statistical tests
\citep[Chapter 3]{Lehmann/Romano:2022}.

\section{Necessity and possibility measures}
\label{sec:confidence_measures}

The notions of a test discussed in the previous sections
allow us to associate measures of confidence with subsets of the parameter space
in view of an outcome.
These are just a different way to package confidence regions.

For a \p-test $(P_\theta\mid\theta\in\Theta)$,
the \emph{\p-necessity measure} of a set $B\subseteq\Theta$
in view of an outcome $\omega\in\Omega$ is defined as
\begin{equation}\label{eq:Box-p}
  \pN(B\mid\omega)
  :=
  \sup_{\theta\notin B}
  P_\theta(\omega).
\end{equation}
Now the Fisher-type disjunction for the true $\theta$ is:
\emph{$\theta\in B$ unless $\omega$ is $\pN(B\mid\omega)$-strange for $\theta$}.
Therefore, we expect $\theta\in B$ for a small $\pN(B\mid\omega)$.
Of course, this disjunction remains true
if we replace ``$\pN(B\mid\omega)$-strange'' by ``$c$-strange'' for any $c\ge\pN(B\mid\omega)$,
but in statements of this kind we usually choose the $c$
that makes them as strong as possible.

Similarly, for an \e-test $(E_\theta\mid\theta\in\Theta)$,
the \emph{\e-necessity measure} of $B\subseteq\Theta$ given $\omega\in\Omega$ is
\begin{equation}\label{eq:Box-e}
  \eN(B\mid\omega)
  :=
  \inf_{\theta\notin B}
  E_\theta(\omega),
\end{equation}
with the analogous interpretation:
\emph{$\theta\in B$ unless $\omega$ is $\eN(B\mid\omega)$-strange for $\theta$}.

If we are interested in a derivative parameter $g(\theta)$, where $g:\Theta\to\Theta_g$,
the \emph{\p-necessity measure} and \emph{\e-necessity measure} of $B\subseteq\Theta_g$
in view of $\omega\in\Omega$ are now defined as
\begin{align}
  \pN^g(B\mid\omega)
  &:=
  \sup_{\theta\in\Theta:g(\theta)\notin B}
  P_\theta(\omega)
  =
  \pN(g^{-1}(B)\mid\omega),
  \label{eq:Box-g-p-used-in-range}\\
  \eN^g(B\mid\omega)
  &:=
  \inf_{\theta\in\Theta:g(\theta)\notin B}
  E_\theta(\omega)
  =
  \eN(g^{-1}(B)\mid\omega),
  \label{eq:Box-g-e}
\end{align}
respectively,
with the same interpretations as before.

Analogously to \eqref{eq:Box-p}--\eqref{eq:Box-g-e}
we can define the \emph{\p-possibility measure} and \emph{\e-possibility measure}
by
\begin{align*}
  \pP(B\mid\omega)
  &:=
  \sup_{\theta\in B}
  P_\theta(\omega)
  =
  \pN(B^{\complement}\mid\omega),\\
  \eP(B\mid\omega)
  &:=
  \inf_{\theta\in B}
  E_\theta(\omega)
  =
  \eN(B^{\complement}\mid\omega),\\
  \pP^g(B\mid\omega)
  &:=
  \pP(g^{-1}(B)\mid\omega)
  =
  \pN^g(B^{\complement}\mid\omega),\\
  \eP^g(B\mid\omega)
  &:=
  \eP(g^{-1}(B)\mid\omega)
  =
  \eN^g(B^{\complement}\mid\omega),
\end{align*}
where $B^{\complement}:=\Theta\setminus B$ is the complement of $B$.
For example, a large value of $\eP(B\mid\omega)$
means that $\theta\in B$ is hardly possible for the true $\theta$
in view of the outcome $\omega$.

\begin{remark}
  Our notation is borrowed from modal logic,
  which has two basic modalities, $\Box$ (necessity) and $\Diamond$ (possibility),
  analogous to the quantifiers $\forall$ and $\exists$, respectively.
  The notions of necessity and possibility measures discussed in this section
  are closely related to the necessity and possibility measures
  of possibility theory \citep{Dubois/Prade:1988}
  (which they include in a wider class of what they call confidence measures),
  and also somewhat related to the belief and plausibility
  functions of the Dempster--Shafer theory \citep{Shafer:1976},
  and to confidence and credibility in conformal prediction
  \citep[(3.66)]{Vovk/etal:2005book}.
  However, unlike their counterparts in those theories,
  our notions just re-express the idea of confidence regions
  without adding new information.
\end{remark}

\subsection*{Necessity measures in the toy example}

In the toy example considered at the end of the previous section
(with the same \p-test and \e-tests),
we can write the \p-necessity measure of a set $B\subseteq\R$ of parameter values
in view of an observation $x\in\R$ as
\begin{equation*}
  \pN(B\mid x)
  =
  2\Phi
  \left(
    -\inf_{\theta\notin B}
    \left|
      x-\theta
    \right|
  \right).
\end{equation*}
According to the definition,
$\pN(B\mid x)$ is determined by the parameter value outside $B$ (assuming the $\inf$ is attained)
that makes the observed $x$ least strange.

The main application of necessity measures in this paper
(described in the following section)
will be ``one-sided'', in that the corresponding confidence regions
will provide only a lower bound
(on the quantity called the number of true discoveries;
equivalently, they provide an upper bound on the number of false discoveries).
If instead of \eqref{eq:E-bar} we use the \e-test \eqref{eq:E-one-sided} with $\delta>0$,
we will have prediction regions in the form of rays pointing left,
and the necessity measure will be
\begin{equation*}
  \eN(B\mid x)
  =
  \exp
  \left(
    \delta
    \left(
      x-\inf B^{\complement}
    \right)
    -
    \delta^2/2
  \right).
\end{equation*}

\section{Controlling the number of false discoveries}
\label{sec:control}

Starting from this section we specialize our setting.
Our sample space $(\Omega,\AAA)$ is still arbitrary,
but now we take $\Theta:=\mathfrak{P}(\Omega)$ as our parameter space
and $Q_\theta:=\theta$ for all $\theta\in\Theta$ as our statistical model;
remember that $\mathfrak{P}(\Omega)$ is the set of all probability measures on $(\Omega,\AAA)$.
Since our statistical model contains all probability measures on $\Omega$,
there is no real loss of generality.

Suppose that we are given $K$ \e-variables $E_1,\dots,E_K$
for testing hypotheses $H_1,\dots,H_K$,
which are our \emph{base hypotheses};
we would like to reject some of them
(in fact, as many of them as possible under a validity constraint).
The realized values of $E_1,\dots,E_K$ are denoted by $e_1,\dots,e_K$,
so that $e_k:=E_k(\omega)$ for the realized outcome $\omega$.

If we do not know anything about the nature of the hypotheses $H_1,\dots,H_K$,
it makes sense to reject a number of them with the largest $e_k$.
But in general, we can consider an arbitrary non-empty \emph{rejection set}
$
  R\subseteq\{1,\dots,K\}
$;
this is the set of base hypotheses, represented by their indices, that the researcher chooses to reject.
\citet{Goeman/Solari:2011local} argue convincingly that in some practically relevant cases
$R$ will not necessarily correspond to the largest $e_k$;
e.g., $R$ may include hypotheses connected by a common theme,
such as all relevant genes related to the gastrointestinal tract
\citep[4.1]{Goeman/Solari:2011local}.

In this section we will find functions $D$ providing a measure of confidence
in the number of true discoveries (to be formally defined momentarily) in the following sense:
a rejection set $R$ contains more than $j$ true discoveries
unless the outcome $\omega$ is $D^R(j)$-strange.
This statement is uniform in $R$ and $j$,
in the sense of the strangeness of outcomes being measured by a fixed \e-test.
Therefore, a large $D^R(j)$ means high confidence
in the number of true discoveries exceeding~$j$.

For each $\theta\in\mathfrak P(\Omega)$, we define
\begin{equation*}
  I_\theta:=\{k\in\{1,\dots,K\}\mid\theta\in H_k\}
\end{equation*}
to be the set of indices of hypotheses containing $\theta$.
If the researcher rejects $H_k$, we refer to this decision as a \emph{discovery}.  
We say that the discovery is \emph{true} if $\theta \notin H_k$,
and it is \emph{false} if $\theta \in H_k$,
where $\theta$ is the true (unknown) probability measure governing the data generation.
For a rejection set $R$, the number of true discoveries is
\begin{equation}\label{eq:g}
  g_R(\theta)
  :=
  \left|
    R\setminus I_\theta
  \right|
  =
  \left|\left\{
    k\in R \mid \theta \notin H_k
  \right\}\right|,
\end{equation}
and the number of false discoveries is
\[
  \left|
    R\cap I_\theta
  \right|
  =
  \left|\left\{
    k\in R \mid \theta \in H_k
  \right\}\right|.
\]
The sum of these two numbers is $\left|R\right|$, the total number of discoveries,
and so controlling the number of false discoveries is the same thing as controlling the number of true discoveries.
Our functions $D^R(j)$ will provide measures of confidence in lower bounds $j+1$ on the number of true discoveries
(equivalently, upper bounds $\left|R\right|-j-1$ on the number of false discoveries).
Researchers are sometimes interested in the proportion of true or false discoveries
$
  \left|R\setminus I_\theta\right| / \left|R\right|
$
or
$
  \left|R\cap I_\theta\right| / \left|R\right|
$,
respectively.
We can also control those with the bounds
$
  (j+1) / \left|R\right|
$
or
$
  (\left|R\right|-j-1) / \left|R\right|
$
(lower for true and upper for false discoveries),
respectively.

\begin{remark}
  The researcher may be interested in parameters $g(\theta)$
  that differ from \eqref{eq:g} more substantially.
  For example, $g(\theta)$ may be the weighted number of true discoveries in $R$
  (e.g., some genes can be more important than other genes).
  Or, for a partition of $R$ into groups
  (one of which can be, e.g., the genes related to the gastrointestinal tract),
  $g(\theta)$ may depend on the number of groups containing true discoveries.
  In this paper we restrict ourselves to the simplest case.
\end{remark}

For \e-confidence bounds, we need an \e-test $(E_\theta)_{\theta\in\mathfrak{P}(\Omega)}$.
For each $k\in I_\theta$,
$E_k$ is an \e-variable for testing $\theta$.
We will obtain $E_\theta$ by merging $(E_k)_{k\in I_\theta}$.
This can be achieved by using \e-merging functions studied in \citet{Vovk/Wang:2021}.
An \emph{\e-merging function} is a Borel function $F:\cup_{n=0}^\infty [0,\infty]^n\to[0,\infty]$
that is increasing in each of its arguments
and maps any finite sequence of \e-variables to an \e-variable:
if $E_1,\dots,E_n$ are \e-variables, $F(E_1,\dots,E_n)$ is required to be an \e-variable as well.
We always set $F:=0$ if the input sequence is empty.
An example (of paramount importance, as discussed earlier) is the arithmetic mean
\begin{equation}\label{eq:AM}
  (e_1,\dots,e_n)
  \mapsto
  \frac1n
  \sum_{i=1}^n
  e_i.
\end{equation}  
An \e-merging function is \emph{symmetric} if it does not depend on the order of its arguments,
like the arithmetic mean.

Let $F$ be a symmetric \e-merging function;
we define for each $\theta\in\Theta$ the \e-variable
\begin{equation}\label{eq:e-test}
  E_\theta:= F(E_k:k\in I_\theta).
\end{equation}
Our main object of interest is $\eP^{g_R}(\{j\}\mid\omega)$ for this \e-test,
which we will abbreviate to $\eP^{g_R}(j\mid\omega)$ dropping the curly braces.

\begin{remark}
  Technically, the choice of $F$ in \eqref{eq:e-test} may even depend on $\theta$,
  but we will ignore this possibility in this paper.
  Moreover, we will be mainly interested in one specific \e-merging function (arithmetic mean).
\end{remark}

Let us replace $\eP^{g_R}(j\mid\omega)$ by a more explicit and easily computable expression.
Set, for a rejection set $R$,
\begin{align}
  \eP^{g_R}(j\mid\omega)
  &=
  \min_{\theta\in\mathfrak{P}(\Omega):g_R(\theta)=j}
  E_\theta
  \notag \\& =
  \min_{\theta\in\mathfrak{P}(\Omega):\left|R\setminus I_\theta\right|=j}
  F(E_k:k\in I_\theta)
  \notag \\ &\ge
  \min_{I\subseteq\{1,\dots,K\}:\left|R\setminus I\right|=j}
  F(E_k:k\in I)
  =:
  D_{\e,F}^R(j),
  \label{eq:D_e_first}
\end{align}
where $\min\emptyset:=\infty$ (as usual),
the argument $\omega$ is implicit after the first ``$=$'',
and the equality $=:$ in \eqref{eq:D_e_first}
signifies $D_{\e,F}^R(j)$ being defined,
with the subscripts (${\e,F}$) dropped if clear from the context.
Intuitively, in \eqref{eq:D_e_first} we go over all $I$
for which there are exactly $j$ true discoveries and evaluate their strangeness;
if all of them are strange,
we are entitled to reject there being exactly $j$ true discoveries.

The values $D^R(j)$ are informative for $j=0,\dots,\lvert R\rvert-1$,
and we will sometimes refer to $D^R(j)$, $j=0,\dots,\lvert R\rvert-1$,
as \emph{discovery \e-vector}.
(Notice that we always have $D^R(\lvert R\rvert)=F(\emptyset)=0$,
and so this value is not informative.)

Let us say (following \citet[Section 1]{Holm:1979})
that $H_1,\dots,H_K$ satisfy the \emph{free combinations condition for $R$}
if the sets $I_\theta$, $\theta\in\Theta$, include all subsets of $R$:
\begin{equation}\label{eq:free}
  \forall S\subseteq R
  \;
  \exists\theta\in\Theta:
  I_{\theta} = S.
\end{equation}
The ``$\ge$'' in \eqref{eq:D_e_first} becomes ``$=$'' under the free combinations condition,
but this condition is not required for the validity of our methods.

\begin{algorithm}[bt]
  \caption{Discovery \e-vector for a given rejection set}
  \label{alg:R}
  \begin{algorithmic}[1]
    \Require
      A symmetric \e-merging function $F$,
      the rejected hypotheses $R\subseteq\{1,\dots,K\}$,
      and an increasing sequence of \e-values $e_1\le\dots\le e_K$.
    \For{$j=0,\dots,\left|R\right|-1$}
      \State let $R_j$ be $R$ without its $j$ largest elements
      \State $D^R(j):=F_{\mathbf{e}}(R_j)$
      \For{$i=1,\dots,\lvert R^{\complement}\rvert$}\label{ln:K-j}
        \State let $R^{\complement}_i$ consist of the $i$ smallest elements of $R^{\complement}$\label{ln:ignore}
        \State $e := F_{\mathbf{e}}(R_j\cup R^{\complement}_i\})$\label{ln:e}
        \If{$e < D^R(j)$}
          \State $D^R(j) := e$
        \EndIf
      \EndFor
    \EndFor
  \end{algorithmic}
\end{algorithm}

An algorithm for computing the discovery vector $D^{R}$ is given as Algorithm~\ref{alg:R};
it is polynomial-time if the underlying \e-merging function $F$, assumed symmetric, is polynomial-time.
It uses the notation $R^{\complement}:=\{1,\dots,K\}\setminus R$ and
\begin{equation}\label{eq:F}
  F_{\mathbf{e}}(I)
  :=
  F(e_i:i\in I),
  \quad
  I\subseteq\{1,\dots,K\},
  \enspace
  I\ne\emptyset,
\end{equation}
where $\mathbf{e}:=(e_1,\dots,e_K)$.
Without loss of generality we assume that the \e-values are sorted in the ascending order,
\begin{equation}\label{eq:sorted}
  e_1\le\dots\le e_K.
\end{equation}

A special and important choice of $F$ is the arithmetic average \eqref{eq:AM}. 
Using this \e-merging function in \eqref{eq:D_e_first},
the \emph{arithmetic-mean discovery \e-vector} is defined as
\begin{equation*}
  \AV^R(j)
  :=
  \min_{I\subseteq\{1,\dots,K\}:\left|R\setminus I\right|=j}
  \frac{1}{\left|I\right|}
  \sum_{i\in I}
  E_i,
  \quad
  j\in\{0,\dots,\left|R\right|-1\}.
\end{equation*}
As we said earlier,
arithmetic averaging is the only useful symmetric \e-merging function
\citep[Proposition~3.1]{Vovk/Wang:2021}.
The vector $\AV^R$ is computed by Algorithm~\ref{alg:R} with
\begin{equation*}
  F_{\mathbf{e}}(I)
  :=
  \frac{1}{\left|I\right|}
  \sum_{i\in I}
  e_i,
  \quad
  I\subseteq\{1,\dots,K\},
  \enspace
  I\ne\emptyset.
\end{equation*}

In general, a discover \e-vector $D_{\e,F}^R(j)$ is not guaranteed to be monotonic in $j$
(not even $\AV^R(j)$ is).
Therefore, we also consider the \emph{regularized} discovery \e-vector
\[
  \bar D_{\e,F}^R(j)
  :=
  \min_{j'\le j}
  D_{\e,F}^R(j').
\]
Regularized discovery \e-vectors automatically satisfy
two other properties of monotonicity.

\begin{proposition}\label{prop:monotone-1}
  For any nonempty sets $R$ and $R'$ in $\{1,\dots,K\}$,
  any $j\in\{0,\dots,\lvert R\rvert-1\}$,
  any $j'\in\{0,\dots,\lvert R'\rvert-1\}$,
  and any \e-merging function $F$:
  \begin{enumerate}[(1)]
  \item\label{it:m_1}
    $\bar D_{\e,F}^R(j')\le\bar D_{\e,F}^R(j)$ if $j\le j'<\lvert R\rvert$;
  \item\label{it:m_2}
    $\bar D_{\e,F}^{R'}(j)\ge\bar D_{\e,F}^R(j)$ if $R\subseteq R'$;
  \item\label{it:m_3}
    $\bar D_{\e,F}^{R'}(j+\left|R'\setminus R\right|)\le\bar D_{\e,F}^R(j)$.
  \end{enumerate}
\end{proposition}

\begin{proof}
  Item \eqref{it:m_1} holds by definition.

  For item \eqref{it:m_2}, we can rewrite the inequality $\bar D_{\e,F}^{R'}(j)\ge\bar D_{\e,F}^R(j)$ as
  \[
    \min_{I:\left|R'\setminus I\right|\le j}
    F(E_i:i\in I)
    \ge
    \min_{I:\left|R\setminus I\right|\le j}
    F(E_i:i\in I),
  \]
  and it suffices to notice that any $I$ satisfying $\left|R'\setminus I\right|\le j$
  satisfies $\left|R\setminus I\right|\le j$.
  (Item \eqref{it:m_2} is also a special case of item \eqref{it:m_3},
  which will be proved independently.)

  Item \eqref{it:m_3} can be rewritten as
  \[
    \min_{I:\left|R'\setminus I\right|\le j+\left|R'\setminus R\right|}
    F(E_i:i\in I)
    \le
    \min_{I:\left|R\setminus I\right| \le j}
    F(E_i:i\in I),
  \]
  which follows from
  \[
    \left|R\setminus I\right| \le j
    \Longrightarrow
    \left|R'\setminus I\right| \le j+\left|R'\setminus R\right|,
  \]
  which in turn follows from the obvious
  \[
    \left|R'\setminus I\right|
    \le
    \left|R\setminus I\right|
    +
    \left|R'\setminus R\right|.
    \qedhere
  \]
\end{proof}

Algorithm~\ref{alg:R} can by adapted to produce the regularized discovery vector
by replacing $\lvert R^{\complement}\rvert$ with $K$ in line~\ref{ln:K-j},
replacing line \ref{ln:e} with $e := F_{\mathbf{e}}(R_j\cup\{1,\dots,i\})$,
and ignoring line \ref{ln:ignore}.

According to our definition of $\bar D$,
\[
  \eN^{g_R}(\{j+1,j+2,\dots\}\mid\omega)
  =
  \eP^{g_R}(\{0,\dots,j\}\mid\omega)
  \ge
  \bar D^R(j).
\]
In agreement with Section~\ref{sec:confidence_measures},
$\bar D^R(j)$ gives a confidence bound on the number of true discoveries:
the rejection set $R$ contains more than $j$ true discoveries
unless the outcome $\omega$ is $\bar D^R(j)$-strange.
Therefore, we can count on there being more than $j$ true discoveries in $R$
for a large observed $\bar D^R(j)$.

Our definitions so far are essentially translations
of Goeman and Solari's [\citeyear[Section~2]{Goeman/Solari:2011local}] definitions
into the language of \e-values.
We will explain the connection in detail in Appendix~\ref{app:GWGS}.
As the procedure for \p-values was first proposed
in \citet[Theorem 6.3]{Genovese/Wasserman:2004},
we refer to it as the \emph{GWGS procedure}.
In Appendix~\ref{app:GWGS} we will also comment on the recent result by \citet{Goeman/etal:2021}
about the GWGS procedure being the only admissible one for controlling true discoveries
(under a property of validity based on \p-values).

\subsection*{Discovery \e-matrices}
\label{sec:matrices}

Next we will discuss a less flexible method
in which we consider a family of rejection sets $R$ that are chosen in an optimal way,
in some sense.
For each $r\in\{1,\dots,K\}$,
the set
\begin{equation}\label{eq:R_r}
  R_r:=\{K-r+1,\dots,K\}
\end{equation}
is the optimal rejection set of size $r$ (assuming \eqref{eq:sorted}),
meaning that $D_{\e,F}^{R_r}\ge D_{\e,F}^R$ for any other set $R\subseteq\{1,\dots,K\}$ of size $r$.
In the terminology of statistical decision theory \citep[Section 1.3]{Wald:1950},
$R_r$ is a complete class of rejection sets.

\begin{algorithm}[bt]
  \caption{Discovery \e-matrix $D$}
  \label{alg:D}
  \begin{algorithmic}[1]
    \Require
      A symmetric \e-merging function $F$
      and an increasing sequence of \e-values $e_1\le\dots\le e_K$.
    \For{$r=1,\dots,K$}
      \For{$j=0,\dots,r-1$}
        \State $S_{r,j}:=\{K-r+1,\dots,K-j\}$
        \State $D_{r,j}:=F_{\mathbf{e}}(S_{r,j})$
        \For{$i=1,\dots,K-r$}\label{ln:K-r}
          \State $e := F_{\mathbf{e}}(S_{r,j}\cup\{1,\dots,i\})$
          \If{$e < D_{r,j}$}
            \State $D_{r,j} := e$
          \EndIf
        \EndFor
      \EndFor
    \EndFor
  \end{algorithmic}
\end{algorithm}

Let us call $D_{r,j}:=D_{\e,F}^{R_r}(j)$ the \emph{discovery \e-matrix}.
Its interpretation is that the outcome $\omega$ is $D_{r,j}$-strange
if there are exactly $j$ true discoveries among the $r$ hypotheses with the largest \e-values.
An algorithm for computing the discovery \e-matrix $D$ is given as Algorithm~\ref{alg:D}.

We are particularly interested in the \emph{arithmetic-mean discovery matrix} $\AM$,
i.e., the discovery \e-matrix
\begin{equation*}
  \AM_{r,j}(e_1,\dots,e_K)
  :=
  \min_{I:\left|R_r\setminus I\right|=j}
  \frac{1}{\left|I\right|}
  \sum_{i\in I}
  e_{i}
  =
  \min_{I:\left|R_r\setminus I\right|\le j}
  \frac{1}{\left|I\right|}
  \sum_{i\in I}
  e_{i}
\end{equation*}
(the last equality will follow from Proposition~\ref{prop:monotone-2}\eqref{it:m_row} below).
In Appendix \ref{app:general} we illustrate discovery \e-matrices with some other choices of the \e-merging function $F$,
which, according to \citet[Proposition 3.1]{Vovk/Wang:2021}, are essentially dominated by the arithmetic mean.

Using the arithmetic mean \e-merging function leads to convenient properties of monotonicity
for the discovery \e-matrix described in the next proposition
(established in Proposition~\ref{prop:monotone-1} above for the regularized version).
The first of these properties will allow us to interpret $D_{r,j}$ in terms of necessity $\Box$,
and the other two properties will help us to visualize $D$ in our plots in the experimental sections.

\begin{proposition}\label{prop:monotone-2}
  For any $R,R'$ of the form \eqref{eq:R_r} (with $r>0$),
  $j\in\{0,\dots,\lvert R\rvert-1\}$,
  $j'\in\{0,\dots,\lvert R'\rvert-1\}$,
  and the arithmetic-mean \e-merging function $F$:
  \begin{enumerate}[(1)]
  \item\label{it:m_row}
    $D_{\e,F}^R(j')\le D_{\e,F}^R(j)$ if $j\le j'<\lvert R\rvert$;
  \item\label{it:m_column}
    $D_{\e,F}^{R'}(j)\ge D_{\e,F}^{R}(j)$ if $R\subseteq R'$;
  \item\label{it:m_diag}
    $D_{\e,F}^{R'}(j+\left|R'\setminus R\right|)\le D_{\e,F}^R(j)$ if $R\subseteq R'$.
  \end{enumerate}
\end{proposition}

\begin{proof}
  The only property of $F$ that we will need is
  \begin{equation}\label{eq:F-mono}
    e\ge\max(\mathbf{e})
    \Longrightarrow
    F(\mathbf{e},e) \ge F(\mathbf{e}),
  \end{equation}
  where $\mathbf{e}\in[0,\infty]^*$ and $e\in[0,\infty]$.
  The arithmetic mean \e-merging function clearly satisfies it.
  It is also clear that \eqref{eq:F-mono} implies
  \begin{equation*}
    \min(\mathbf{e}')\ge\max(\mathbf{e})
    \Longrightarrow
    F(\mathbf{e},\mathbf{e}') \ge F(\mathbf{e}),
  \end{equation*}
  for any $\mathbf{e},\mathbf{e}'\in[0,\infty]^*$.

  Item \eqref{it:m_row} can be rewritten as
  \begin{equation}\label{eq:m_row}
    \min_{I':\left|R\setminus I'\right|=j'}
    F(E_i:i\in I')
    \le
    \min_{I:\left|R\setminus I\right|=j}
    F(E_i:i\in I).
  \end{equation}
  Let $I=R_j\cup R^{\complement}_i$ be a set
  where the $\min$ on the right-hand side of \eqref{eq:m_row} is attained
  as Algorithm~\ref{alg:R} is run.
  It suffices to consider $I':=R_{j'}\cup R^{\complement}_i$.

  Items \eqref{it:m_column} and \eqref{it:m_diag} follow
  from item \eqref{it:m_row} of this proposition
  in combination with items \eqref{it:m_2} and \eqref{it:m_3} of Proposition~\ref{prop:monotone-1}.

  Properties~\eqref{it:m_row}--\eqref{it:m_diag} are not independent;
  namely, \eqref{it:m_row} follows from \eqref{it:m_column} and \eqref{it:m_diag}.
  Indeed, given $R$, $j$, and $j'$ as in \eqref{it:m_row}
  and choosing any $R'\supseteq R$ satisfying $\left|R'\setminus R\right|=j'-j$,
  we can deduce from \eqref{it:m_column} and \eqref{it:m_diag}:
  \[
    D^R(j')
    =
    D^R(j+\left|R'\setminus R\right|)
    \le
    D^{R'}(j+\left|R'\setminus R\right|)
    \le
    D^R(j).
    \qedhere
  \]
\end{proof}

Proposition~\ref{prop:monotone-2}\eqref{it:m_row} gives us an interpretation of the discovery \e-matrix
in terms of the Fisher-type disjunction:
there are more than $j$ true discoveries among the $r$ hypotheses with the largest \e-values
unless the outcome $\omega$ is $D_{r,j}$-strange.
In symbols,
\[
  \eN^{g_{R_r}}(\{j+1,j+2,\dots\}\mid\omega)
  \ge
  D_{r,j}(\omega).
\]
Proposition~\ref{prop:monotone-2} also shows that discovery matrices $D_{r,j}$ are monotonic functions along the rows $r=\const$,
along the columns $j=\const$, and along the diagonals $r-j=\const$.
Notice that the monotonicity along the rows follows immediately
from the monotonicity along the columns and the monotonicity along the diagonals.

\begin{remark}\label{rem:independent}
  It is true that arithmetic averaging is the only useful symmetric \e-merging function
  when no assumptions are made about the dependence structure of the base \e-values.
  On the other hand, if the base \e-values are supposed to be independent,
  it is clear that the product of \e-variables is always an \e-variable.
  Therefore, when defining the \e-test \eqref{eq:e-test},
  we have plenty of alternatives to the arithmetic mean in the role of $F$;
  $F$ can be the product, or a combination of the arithmetic mean and the product.
  The product does not work well for problems of the type considered in this paper,
  since small base \e-values then have disproportionate effect.
  However, the combination
  \begin{equation*}
    E_\theta
    :=
    \sum_{\{i,j\}\subseteq I_\theta:i\ne j}
    E_i E_j
    /
    \binom{\left|I_\theta\right|}{2}
  \end{equation*}
  (with $\{i,j\}$ ranging over the 2-element subsets of $I_\theta$)
  of the product and arithmetic averaging
  gives excellent results,
  much better than what we can get without the assumption of independence.
  See \citet{Vovk/Wang:eWP4} for details.
\end{remark}

\section{Simulation studies}
\label{sec:simulation}

In our simulation studies we will visualize the arithmetic-mean discovery matrix in some simple cases
and compare Algorithm~\ref{alg:D} with a method based on \p-values.
Our setting will be similar to that of \citet[Section~5]{Vovk/Wang:2021},
where family-wise validity is studied.

The observations are generated from the Gaussian model $N(\mu,1)$.
The null hypotheses are $N(0,1)$ and the alternatives are $N(\delta,1)$,
where we take $\delta:=-3$ throughout the section.
We generate $K/2$ observations from $N(\delta,1)$ (the alternative distribution)
and then $K/2$ observations from $N(0,1)$ (the null distribution),
where $K$ (the overall number of hypotheses) is an even number.

\begin{figure}
  \begin{center}
    \includegraphics[width=0.6\textwidth]{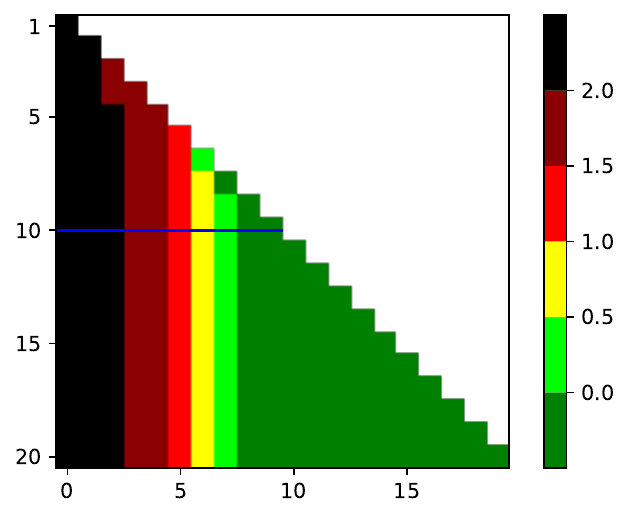}
  \end{center}
  \caption{The arithmetic-mean discovery matrix for 10 false
    and 10 true null hypotheses,
    as described in text.
    The colour map on the right gives Jeffreys's thresholds,
    the boundaries between different colours in most of our plots,
    on the decimal log scale.
   Row 10 is highlighted in blue.}
  \label{fig:small}
\end{figure}

In this paper,
we colour-code the entries (\e-values) of discovery \e-matrices
according to Jeffreys's rule of thumb discussed in Section \ref{sec:testing}.
The full colour map is shown on the right of Figure~\ref{fig:small}
with the thresholds between different colours given in terms of the decimal logarithm of \e-values.
The most interesting parts of our plots of discovery \e-matrices are those in yellow and red;
green and dark green parts carry little or no evidence and so are useless for us,
and dark red and black parts carry so much evidence that they are rare in a wide range of practical applications
(cf.\ Section~\ref{sec:empirical}).
In all our discussions below we will ignore the boundaries between the green and dark green parts.

Figure~\ref{fig:small} shows the arithmetic-mean discovery matrix that Algorithm~\ref{alg:D} gives
for $K=20$: we generate $10$ observations from $N(\delta,1)$ and then $10$ from $N(0,1)$.
The base \e-values are the likelihood ratios
\begin{equation}\label{eq:E}
  E(x)
  :=
  \frac{\dd N(\delta,1)}{\dd N(0,1)}(x)
  =
  \exp(\delta x - \delta^2/2)
\end{equation}
(cf.\ \eqref{eq:E_mu})
of the alternative to the null density,
where $x\sim N(\mu,1)$ is the corresponding observation.
For example, row 10 (highlighted in blue) of the matrix in Figure~\ref{fig:small} shows
that there is decisive evidence that the number of true discoveries
among the 10 hypotheses with the largest \e-values is at least 3.
Similarly, there is very strong evidence that the number of true discoveries is at least 5,
there is strong evidence that the number of true discoveries is at least 6, etc.

In this and following sections we will see many representations of discovery \e-matrices
resembling Figure~\ref{fig:small}.
They are also convenient representations of confidence regions for the numbers of true discoveries
at various significance levels.
For example, the non-black part of each such figure is formed by the extremely strong confidence regions in each row.
In the case of Figure~\ref{fig:small},
row 10 shows that the extremely strong confidence region for the number of true discoveries
among the 10 hypotheses with the largest \e-values is $\{3,\dots,10\}$.
The red/yellow/green part (not including dark red) is formed by the very strong confidence regions,
so that the very strong confidence region for the number of true discoveries
among the 10 hypotheses with the largest \e-values is $\{5,\dots,10\}$.
Similarly, the yellow/green and green parts are formed by strong and substantial confidence regions,
respectively.

\begin{figure}
  \begin{center}
    \includegraphics[width=0.32\textwidth]{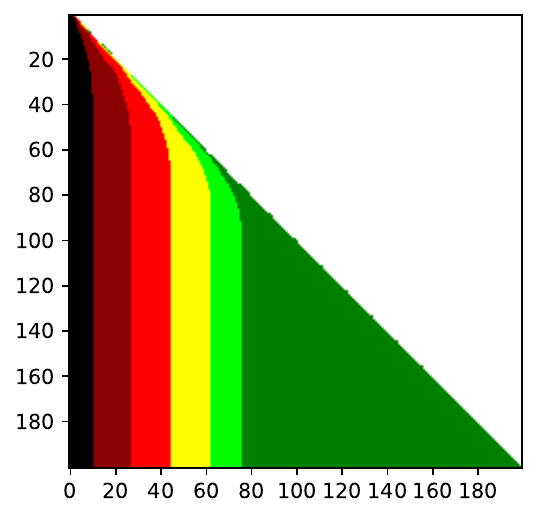}
    \includegraphics[width=0.32\textwidth]{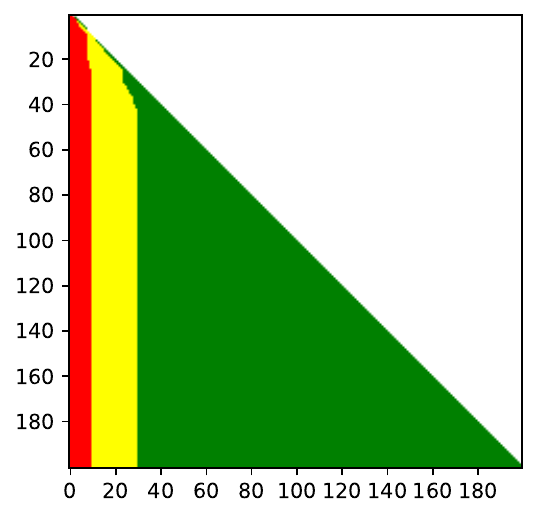}
    \includegraphics[width=0.32\textwidth]{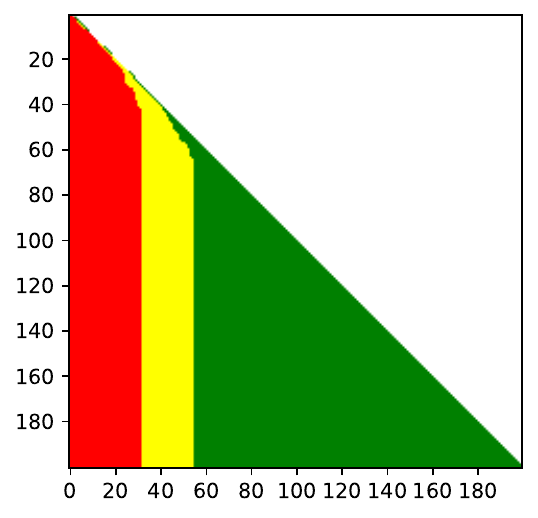}\\
    \includegraphics[width=0.32\textwidth]{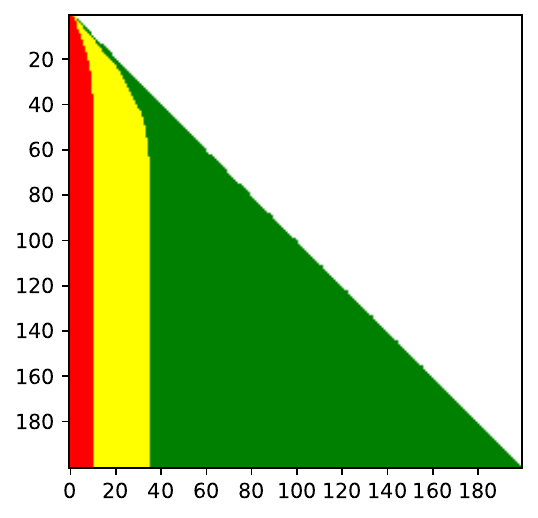}
    \includegraphics[width=0.32\textwidth]{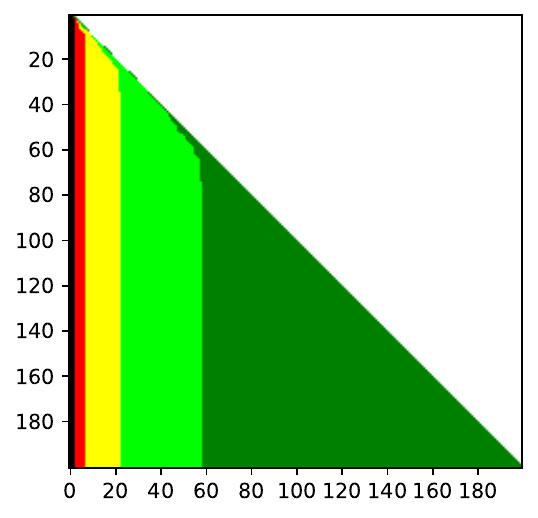}
    \includegraphics[width=0.32\textwidth]{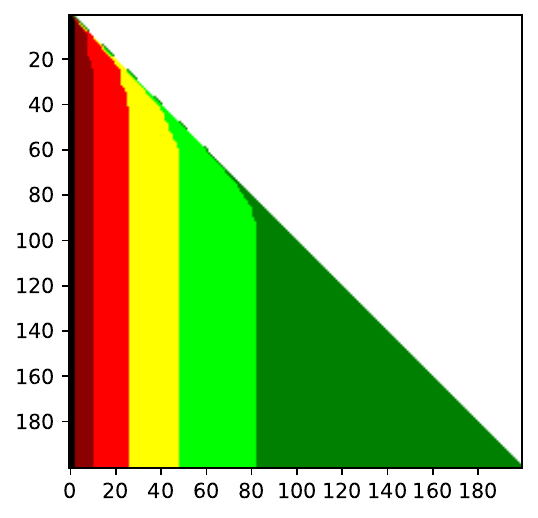}
  \end{center}
  \caption{Upper left panel:
    the arithmetic-mean discovery matrix for 100 false and 100 true null hypotheses.
    Upper middle panel:
    the GWGS discovery \p-matrix in the same situation
    for Fisher's thresholds $1\%$ and $5\%$
    (with values below $1\%$ shown in red and between $1\%$ and $5\%$ in yellow)
    under arbitrary dependence.
    Upper right panel:
    as the upper middle panel but assuming independence.
    \newline
    Lower left panel:
    the \e-to-\p\ calibrated arithmetic-mean discovery matrix
    in the upper left panel
    using Fisher's thresholds.
    Lower middle panel:
    the VS-transformed GWGS discovery \p-matrix
    in the upper middle panel of Figure~\ref{fig:big} (under arbitrary dependence)
    using Jeffreys's thresholds.
    Lower right panel:
    the VS-transformed GWGS discovery \p-matrix
    in the upper right panel
    (under independence)
    using Jeffreys's thresholds.}
  \label{fig:big}
\end{figure}

The upper left panel of Figure~\ref{fig:big} is the counterpart of Figure~\ref{fig:small} for a larger number of hypotheses, $K=200$:
we generate $100$ observations from $N(\delta,1)$ and then $100$ from $N(0,1)$.
We will refer to this set of observations as the \emph{simulation data}.

\begin{table}
  \caption{The values $D_{r,j}$ of the discovery matrix shown in Figure~\ref{fig:big}
    for several rows $r$ and columns $j$.\label{tab:D}}
  \begin{center}{\footnotesize\!\!\!\!\!
    \begin{tabular}{c|cc|ccc|cccc|cc}
      \diagbox{$r$}{$j$} & 10 & 11 & 26 & 27 & 28 & 39 & 40 & 41 & 42 & 45 & 46 \\
      \hline
      31 & 95.6 & 87.5 & 14.1 & \textbf{10.8} & \textbf{7.63} & & & & & & \\
      32 & 96.5 & 88.5 & 16.0 & \textbf{12.7} & \textbf{9.63} & & & & & & \\
      \hline
      50 & \textbf{103} & \textbf{96.0} & \textbf{33.2} & \textbf{30.6} & 28.2 & \textbf{11.2} & \textbf{9.94} & 8.73 & 7.52 & \textbf{4.04} & \textbf{3.11} \\
      51 & 103 & 96.0 & 33.4 & 30.9 & 28.5 & 11.6 & \textbf{10.4} & \textbf{9.21} & 8.01 & 4.61 & 3.69 \\
      52 & 103 & 96.0 & 33.6 & 31.1 & 28.7 & 12.0 & \textbf{10.8} & \textbf{9.61} & 8.43 & 5.10 & 4.19
    \end{tabular}}
  \end{center}
\end{table}

In practice, a discovery \e-matrix, such as that shown in the upper left panel of Figure~\ref{fig:big},
can be used in different ways, for example:
\begin{itemize}
\item
  The researcher may have budget for a limited number of follow-up studies of the hypotheses.
  For example, if in the situation of that panel
  her budget is 50 hypotheses,
  she just concentrates on row 50 (studying the 50 hypotheses $H_k$ with the largest \e-values).
  For the first 11 entries (namely, those with indices 0 to 10 inclusive) in this row the \e-value exceeds $100$,
  and so she has decisive evidence that there are at least $11$ true discoveries among those 50 hypotheses.
  Similarly,
  \begin{itemize}
  \item
    she has very strong evidence that there are at least $27$ true discoveries,
  \item
    she has strong evidence that there are at least $40$ true discoveries,
  \item
    she has substantial evidence that there are at least $46$ true discoveries.
  \end{itemize}
  For the relevant \e-values, see the bold entries in the row $r=50$ of Table~\ref{tab:D}.
  In terms of confidence regions, we can say, e.g.,
  that our method gives the substantial \e-confidence region $\{46,47,\dots\}$,
  so that 46 may be called the substantial lower \e-confidence bound on the number of true discoveries
  among the 50 hypotheses.
\item
  The researcher might have some idea of what proportion of false discoveries she is willing to tolerate
  (in the spirit of choosing the false discovery rate \emph{a priori}
  \citep{Benjamini/Hochberg:1995}).
  For example, if she is willing to tolerate $10\%$ of false discoveries
  and willing to use Jeffreys's standard (\e-value greater than 10) of strong evidence,
  she should concentrate on row 31 (i.e., study the 31 hypotheses with the largest \e-values),
  which is the lowest row with at most $10\%$ of entries below $10$.
  See the bold entries in the rows $r\in\{31,32\}$ of Table~\ref{tab:D}
  (we have strong evidence that
  there are at most $3/31\approx9.7\%$ of false discoveries in row 31
  and at most $4/32=12.5\%$ of false discoveries in row 32).
\item
  Alternatively, the researcher might have some idea of how many false discoveries she is willing to tolerate
  (in the spirit of $k$-FWER \citep{Romano/Wolf:2007}).
  If she is willing to tolerate at most $10$ false discoveries
  and still willing to use Jeffreys's standard of strong evidence,
  she should concentrate on row 51,
  which is the lowest row with at most $10$ entries
  (in fact, exactly $10$ entries)
  below $10$.
  See the bold entries in the rows $r\in\{51,52\}$ of Table~\ref{tab:D}
  (we have strong evidence that there are at most $10$ false discoveries in row 51
  and at most $11$ false discoveries in row 52).
\end{itemize}
Of course, the researcher may know her hypotheses and relations between them very well,
and after looking at the discovery \e-matrix she may come up with her own rejection set $R$,
as discussed in Section~\ref{sec:control}.
In this case she should also use Algorithm~\ref{alg:R}.

All discovery \e-matrices in this and following sections,
as noticed above (after Proposition~\ref{prop:monotone-2}),
satisfy three (non-independent) properties of monotonicity. 
Namely, the entries $\AM_{r,j}$ are decreasing in $j$ (for a fixed $r$),
increasing in $r$ (for a fixed $j$),
and decreasing in $r$ and $j$ along the lines parallel to the main diagonal.

\subsection*{Comparisons}

This paper concentrates on multiple hypothesis testing using \e-values,
but in scientific practice \p-values are more popular,
despite recent criticism.
In this subsection we will report results of our simulation studies in terms of \p-values
and compare them to our results, as best we can
in view of the difficulties discussed in Section~\ref{sec:comparison}.

For comparison with methods based on \p-values,
we use the GWGS procedure applied to standard procedures for combining \p-values
and to the same nested rejection sets
(initial subsets of $\{1,\dots,200\}$ assuming the \p-values are given in ascending order).
These procedures admit computationally efficient shortcuts \citep[Theorem~1]{Goeman/etal:2019Biometrika}
and are implemented in the \textsf{R} package \texttt{hommel} \citep{Goeman/etal:2019R}.
As the base \p-values we take
$
  P(x)
  :=
  \Phi(x)
$,
where $\Phi$ is, as before, the standard Gaussian distribution function;
these are the \p-values found using the most powerful test given by the Neyman--Pearson lemma.
The GWGS procedure can be interpreted as producing an analogue of a discovery \e-matrix,
which we call a \emph{discovery \p-matrix},
with \e-values replaced by \p-values.
For details, see Appendix~\ref{app:GWGS}.
In particular, a version of the notion of a discovery \p-vector
was introduced in \citet[Section~3]{Goeman/Solari:2011rejoinder}.

The package \texttt{hommel} has an option (\texttt{simes})
that controls the choice of the procedure for combining \p-values,
and the resulting discovery \p-matrix is valid either under arbitrary dependence,
in which case Hommel's \citeyr{Hommel:1986} procedure is used for combining \p-values,
or under certain assumptions on the dependence structure for the input \p-values,
in which case Simes's \citeyr{Simes:1986} procedure is used.
In particular, Simes's procedure is valid under the assumption of independence;
it is also valid under relaxations of independence such as positive dependence \citep{Sarkar:2011},
but not under arbitrary dependence.
For brevity we will talk about \p-values that are either arbitrarily dependent or independent,
but it should be remembered that the assumption of independence may be relaxed.
We never make such assumptions about base \e-values (but cf.\ Remark~\ref{rem:independent}).

The upper middle panel of Figure~\ref{fig:big} shows
the discovery \p-matrix found using \texttt{hommel} applied to the simulation data
under arbitrary dependence.
The upper right panel of Figure~\ref{fig:big} is analogous but assumes independent base \p-values.
Both panels use Fisher's thresholds $1\%$ and $5\%$;
the values below $1\%$ are shown in red,
between $1\%$ and $5\%$ in yellow, and above $5\%$ in green
(so that red means ``highly significant'' and yellow means ``significant but not highly significant'').
According to Jeffreys as quoted in Section~\ref{sec:comparison} (p.~\pageref{p:Jeffreys}),
the red and yellow areas are somewhat comparable between \p-values and \e-values,
but we can draw some conclusions even without such cross-comparisons.

Remember that our method does not require any assumptions about the dependence structure of the \e-values.
It is true that our simulated data are independent,
but this information is typically unavailable,
and the performance of methods that do not depend on independence or similar assumptions is still interesting.
Comparing the upper left and upper middle panels of Figure~\ref{fig:big},
we can see that our method produces better confidence bounds
if we are willing to use Jeffreys's informal correspondence between \e-values and \p-values.
The upper left panel is even better, in this sense, than the upper right panel,
which makes an assumption on the dependence structure of the base \p-values.

The three lower panels of Figure~\ref{fig:big} are the transformed versions of the corresponding upper panels.
In the lower left panel, we transform the arithmetic-mean discovery matrix (upper left panel)
by applying the canonical \e-to-\p\ calibrator $e\mapsto1/e$.
In the other two lower panels, we transform the corresponding upper panels
by applying the VS transformation \eqref{eq:VS}.
Therefore, the lower left panel contains valid \p-values, whereas the other two lower panels
contain upper bounds on \e-values.

It is interesting that even after the crude step of \e-to-\p\ calibration
(remember the woeful round-trip efficiency illustrated by \eqref{eq:round_trip}),
the lower left panel of Figure~\ref{fig:big} still looks slightly better than the upper middle panel.
In this comparison there is no uncertainty in the choice of the \e-to-\p\ calibrator,
since \eqref{eq:e-to-p_calibrator} is the only reasonable one
(namely, it dominates any other \e-to-\p\ calibrator).
And even the optimistic VS transformation (the lower middle panel)
looks much worse than the arithmetic-mean discovery matrix in the upper left panel.
We can see, even without using Jeffreys's informal correspondence,
that the method based on \e-values produces better results in this case,
despite the crude calibration steps.

Not surprisingly,
assuming independence makes direct treatment of \p-values more efficient:
compare the upper right panel and the lower left panel.
What is more surprising is that, even assuming independence and using the optimistic VS transformation,
the lower right panel still look worse than the upper left panel.

\begin{table}
  \caption{The Benjamini--Hochberg and Benjamini--Yekutieli procedures
    applied to the simulation data for FDR (false discovery rate) $5\%$ and $1\%$.\label{tab:sim_BH}}
  \begin{center}
    \begin{tabular}{lrr}
      assumption & $5\%$ & $1\%$ \\
      \hline
      independence & 87 & 61 \\
      arbitrary dependence & 55 & 28
    \end{tabular}
  \end{center}
\end{table}

Table~\ref{tab:sim_BH} gives the numbers of null hypotheses rejected by the Benjamini--Hochberg procedure
\citep{Benjamini/Hochberg:1995}
and its version for arbitrary dependence
\citep{Benjamini/Yekutieli:2001}.
(The relationship between the Benjamini--Yekutieli procedure and the \texttt{hommel} package without the \texttt{simes} option
is the same as that between the Benjamini--Hochberg procedure and the \texttt{hommel} package with the \texttt{simes} option;
this is made explicit in \citet[Section 5]{Goeman/etal:2019Biometrika}.)
Here the results are more difficult to interpret,
since the kind of guarantees provided by those procedures is so different
from the guarantees provided by the GWGS methods.
Roughly, we get comparable results between the row ``arbitrary dependence'' in Table~\ref{tab:sim_BH}
and the confidence bounds for the number of true discoveries in, say, row 50
of the discovery matrix in the upper left of Figure~\ref{fig:big}:
as discussed earlier (cf.\ Table~\ref{tab:D}),
the strong and substantial confidence bounds are 40 and 46, respectively.
(Remember that, following Jeffreys, ``strong'' refers to the threshold of $10$ for \e-values
and regarded as roughly corresponding to ``highly significant'',
and ``substantial'' refers to the threshold of $10^{1/2}$ for \e-values
and regarded as roughly corresponding to ``significant''.)

\section{Empirical studies}
\label{sec:empirical}

In this section we will demonstrate how the methods of this paper can, in principle, be used in practice.
It is important that we will make no assumptions of independence.

We will use the classical dataset first described in \citet{Hedenfalk/etal:2001}
and then carefully studied in \citet{Storey/Tibshirani:2003}.
Essentially, we will adapt Storey and Tibshirani's analysis to using \e-values in place of \p-values
(see the end of the section for a discussion of differences).
In our experiments we use the version of the dataset made available as part of the \textsf{R} package \texttt{qvalue}
\citep{Storey/etal:2019}.

The main content of the dataset is the expression levels of $3226$ genes
in 15 samples of tissues.
Seven samples are coming from carriers of mutations in the BRCA1 gene,
and the remaining eight from carriers of mutations in the BRCA2 gene.
We will say that each sample is labelled with its BRCA status:
seven are labelled BRCA1, and eight are labelled BRCA2.
The core of the dataset is the $3226\times15$ matrix
of gene expressions in the samples;
all entries are positive numbers.
Following \citet[Appendix, Remark C]{Storey/Tibshirani:2003},
we remove all rows containing at least one entry exceeding $20$,
which leaves us with a $3170\times15$ data matrix.
Each row of the data matrix corresponds to a gene
and each column to a sample,
and the entry in row $k$ and column $j$ is the expression level of gene $k$ in sample $j$.
For each gene we are interested in the scientific hypothesis
that the gene expression does not depend on the BRCA status of the sample.

Storey and Tibshirani's version of the dataset also contains some further information,
such as the \p-value for each gene.
For further information about this dataset,
which we will refer to as the \emph{BRCA dataset},
see, e.g., \citet[Section 5]{Storey/etal:2007} and \citet[5.2]{Guindani/etal:2009}.

For this dataset methods ensuring family-wise validity do not work well.
For example, the ten smallest \p-values in Storey and Tibshirani's list
multiplied by the number of genes 3170 are
\begin{equation}\label{eq:10_p}
  0, \;
  0.040, \;
  0.050, \;
  0.060, \;
  0.060, \;
  0.100, \;
  0.140, \;
  0.160, \;
  0.240, \;
  0.240
\end{equation}
and so the Bonferroni correction leads to only three statistically significant \p-values.
Moreover, one of the \p-values is exactly zero,
and so cannot be a valid \p-value
(for details, see p.~\pageref{p:zero_p-value}).
Hedenfalk et al.\ conclude that 9--11 genes are differentially expressed.
Storey and Tibshirani's informal analysis suggests that many more,
at least 33\%,
of the examined genes are differentially expressed.
However, their informal analysis assumes what they call ``weak independence'':
they rely on the law of large numbers when inspecting histograms of \p-values,
assuming that the probabilities of the \p-values lying in various ranges
will manifest themselves as empirical frequencies seen in the histograms.
Their formal analysis is asymptotic and also assumes weak independence:
see their Appendix, Remark~D.

We formalize the scientific theory of interest
as the following statistical hypothesis about each gene $k$:
given the multiset of expression levels in row $k$ of the data matrix,
each ordering of the row has the same probability.
In our current context,
a \emph{nonconformity measure} is a measurable function of two multisets;
we will use it by applying, for a given gene,
to the multiset (of size 7) of the expression levels for the samples labelled BRCA1
and the multiset (of size 8) of the expression levels for the samples labelled BRCA2;
the resulting value will be called the \emph{nonconformity score}.
For computing base \e-values, we use the formula
\begin{equation}\label{eq:permutation_e}
  e_k
  :=
  \frac{T_k}
  {
    \frac{1}{B + 1}
    \left(
      \sum_{b=1}^B
      T_k^{0b}
      +
      T_k
    \right)
  },
  \quad
  k=1,\dots,3170,
\end{equation}
where $T_k$ is the nonconformity score computed from the $k$th row of the data matrix
with the true labels (BRCA1 or BRCA2) for each sample,
$T_k^{0b}$ is the nonconformity score computed from the same row with randomly permuted labels,
and $B$ is the number of permutations.
In our experiments, the case $0/0$ of the right-hand side of \eqref{eq:permutation_e} never occurs.
We will call \eqref{eq:permutation_e} the \emph{Monte Carlo \e-value}.
We are justified in calling it an \e-value since,
under the null hypothesis,
the expected value of the right-hand side of \eqref{eq:permutation_e} is 1
if we set $0/0:=1$.
(Moreover, the conditional expectation of the right-hand side of \eqref{eq:permutation_e} is 1
given the multiset of nonconformity scores $\{T_k^{01},\dots,T_k^{0B},T_k\}$,
the expectation being over all choices of the position of the true nonconformity score in the multiset.)

Our nonconformity measure will be defined in terms of the \textit{t}-statistic.
Let $x_{k j}$ be the base two logarithm of the value in row $k$ and column $j$ of the data matrix
(although the base does not matter in our empirical studies).
The two-sample \textit{t}-statistic for the $k$th gene is
\begin{equation}\label{eq:t}
  t_k
  :=
  \frac{\bar x_{k2} - \bar x_{k1}}{\sqrt{s^2_{k1} / n_1 + s^2_{k2} / n_2}},
\end{equation}
where $n_1=7$ is the number of BRCA1 columns,
$n_2=8$ is the number of BRCA2 columns,
and
\[
  \bar x_{k1}
  :=
  \frac{1}{n_1}
  \sum_{j\in\BRCAI}
  x_{k j},
  \quad
  s^2_{k1}
  :=
  \frac{1}{n_1-1}
  \sum_{j\in\BRCAI}
  \left(
    x_{k j} - \bar x_{k1}
  \right)^2
\]
are the sample mean and variance for the BRCA1 entries,
with the analogous expressions for BRCA2.
The variances of the two groups (BRCA1 and BRCA2) are not assumed to be equal
(following \citet{Storey/Tibshirani:2003}),
but using equal-variance two-sample \textit{t}-statistics would lead to similar results.

We define the nonconformity score as $T_k:=f(t_k)$ for some function $f$ of the \textit{t}-statistic $t_k$
(see \eqref{eq:t}).
A natural nonconformity score is $\left|t_k\right|$,
but we generalize it to $\left|t_k\right|^d$ for some $d>0$.
This choice of $f$ is motivated by the Bayesian two-sample \textit{t}-test
widely discussed in recent literature starting from \citet{Gonen/etal:2005}
and briefly reviewed in \citet[Section~3]{Gonen/etal:2019}.
A standard expression for the Bayes factor produced by such a test
via the \textit{t}-statistic $t$ is
\begin{equation}\label{eq:expression}
  f(t)
  :=
  c
  \left(
    1+a t^2
  \right)^{d/2}
\end{equation}
for positive constants $a$, $c$, and $d$ involving the number of degrees of freedom and effective sample size;
see, e.g., \citet[(14)]{Wang/Liu:2016} and \citet[(1)]{Rouder/etal:2009};
the form \eqref{eq:expression} goes back to Jeffreys \citep[(12)]{Ly/etal:2016}.
However, different constants are used in different papers.
We set, without loss of generality, $c:=1$,
since $c$ cancels out when using \eqref{eq:permutation_e}.
We further simplify \eqref{eq:expression} by ignoring the ``$1+{}$'';
this makes $a$ and any constant factors in the definition of the \textit{t}-statistic $t$
(there is a non-trivial factor under the assumption of equal variances for the two groups)
irrelevant,
as they also cancel out when applying \eqref{eq:permutation_e}.
Of course, this step does not affect the validity of our methods.

\begin{figure}
  \begin{center}
    \includegraphics[width=0.48\textwidth]{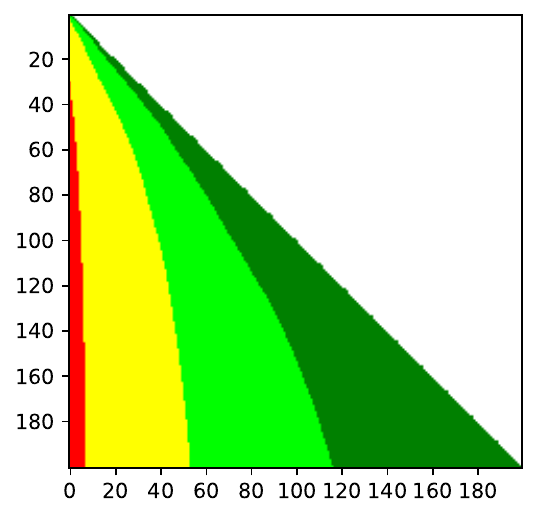}
    \includegraphics[width=0.48\textwidth]{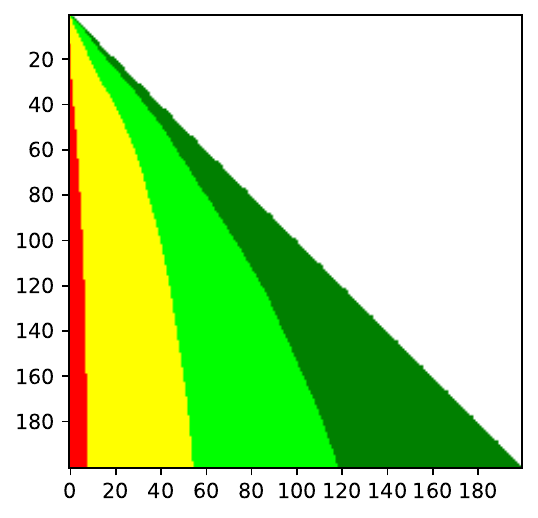}
  \end{center}
  \caption{Left panel:
    the top-left $200\times200$ corner of the arithmetic-mean discovery matrix
    for the \texttt{BRCA} dataset for $B:=10000$,
    using Jeffreys's thresholds.
    Right panel: its version (based on \eqref{eq:permutation_e_simple})
    that is only approximately valid.}
  \label{fig:BRCA_10000}
\end{figure}

The left panel of Figure~\ref{fig:BRCA_10000}
gives a key part of the arithmetic-mean discovery matrix for the BRCA dataset
with $f(t):=\left|t\right|^d$ for $d:=10$, with $B:=10000$,
and with base Monte Carlo \e-values \eqref{eq:permutation_e}.
We can see that there is strong evidence that the number of differentially expressed genes
is at least as large as Hedenfalk at al.'s number.
If we settle for substantial evidence, the number is much larger.
Arguably, it is not as large as in Storey and Tibshirani's study,
but we are not using any exchangeability or independence assumptions.

\begin{table}
  \caption{Summaries of the last row of the arithmetic-mean discovery matrix
    for different values of $d$.
    Column ``strong'' contains the number of entries that are greater than 10
    (all of them are below $10^{3/2}$, and so they provide strong evidence, i.e., red in our pictures).
    Column ``at least substantial'' contains the number of entries that are greater than $10^{1/2}$
    (providing at least substantial evidence).\label{tab:emp_d}}
  \begin{center}
    \begin{tabular}{ccc}
      $d$ & strong & at least substantial \\
      \hline
      4 & 0 & 62 \\
      6 & 0 & 82 \\
      8 & 4 & 70 \\
      10 & 7 & 56 \\
      12 & 8 & 46 \\
      20 & 9 & 29 \\
      50 & 8 & 17 \\
      100 & 7 & 14
    \end{tabular}
  \end{center}
\end{table}

Dependence on the initial state of the random numbers generator
(always set to 1 in our experiments) is fairly significant
but does not affect our conclusions.
Dependence on the value of $d$ is also significant;
the values below 10 tend to lead to higher numbers of true discoveries for Jeffrey's standard of substantial evidence,
and the values above 10 to higher numbers of true discoveries for strong evidence
(up to a limit; see Table~\ref{tab:emp_d}).
The literature on the Bayesian two-sample \textit{t}-test quoted above
seems to suggest that $d$ should have the same order of magnitude as the number of degrees of freedom.

\subsection*{Comparisons}

We start by comparing our methodology with that of \citet{Storey/Tibshirani:2003},
which was our main source of data and ideas in this section.
The main differences are:
\begin{itemize}
\item
  Storey and Tibshirani use \p-values whereas we use \e-values.
\item
  Storey and Tibshirani implicitly assume that the genes are exchangeable
  under the null hypothesis.
\item
  Moreover, Storey and Tibshirani assume that the \p-values are weakly independent.
\end{itemize}

\label{p:zero_p-value}
Strictly speaking, Storey and Tibshirani's method does not produce valid \p-values,
even under their null hypothesis implicitly involving gene exchangeability.
This can be seen from their formula for computing the \p-values,
\begin{equation}\label{eq:ST}
  p_k
  :=
  \frac
  {
    \sum_{b=1}^B
    \left|\{j : \left|t_j^{0b}\right|\ge\left|t_k\right|, j=1,\dots,3170\}\right|
  }
  {3170\cdot B}
\end{equation}
(the last displayed equation in their Appendix, Remark C),
where $t_k$ is the \textit{t}-statistic for gene $k$
and $t_j^{0b}$ is the \textit{t}-statistic for gene $j$ with the labels BRCA1 and BRCA2 randomly permuted
(for the $b$th random permutation, $b=1,\dots,B$ and $B:=100$).
The numerator of \eqref{eq:ST} can well be zero
(and it is in one case: see \eqref{eq:10_p}).

To turn the expression \eqref{eq:ST} into a valid \p-value
(under the null hypothesis of gene exchangeability and label uninformativeness),
it suffices to add $1$ to the numerator and denominator of \eqref{eq:ST};
cf.\ \citet[(17.7)]{Lehmann/Romano:2022}, \citet{Hemerik/Goeman:2018},
and the method of conformal prediction \citep{Vovk/etal:2005book}.
Namely,
\begin{equation}\label{eq:p-conformal}
  \frac
  {
    \sum_{b=1}^B
    \left|\{j : \left|t_j^{0b}\right|\ge\left|t_k\right|, j=1,\dots,3170\}\right|
    +
    1
  }
  {3170\cdot B + 1}
\end{equation}
is a valid \p-value.
The intuition behind the expression~\eqref{eq:p-conformal}
is that, to see how well $t_k$ conforms to the multiset of size $3170\cdot B$
consisting of $t_j^{0b}$,
we add $t_k$ to the multiset before computing the rank \p-value.
Since we are comparing the \textit{t}-statistic for gene $k$ with \textit{t}-statistics for other genes 
in~\eqref{eq:ST} and~\eqref{eq:p-conformal},
we are implicitly assuming gene exchangeability.

Under gene exchangeability, for computing base \e-values, we can use the formula
\begin{equation*}
  e_k
  :=
  \frac{T_k}
  {
    \frac{1}{3170\cdot B + 1}
    \left(
      \sum_{b=1}^B
      \sum_{j=1}^{3170}
      T_j^{0b}
      +
      T_k
    \right)
  },
\end{equation*}
in analogy with \eqref{eq:p-conformal}.
This gives an \e-variable under the assumption that the labels are uninformative
and the genes are exchangeable.

To avoid the assumption of gene exchangeability,
we use the expression \eqref{eq:permutation_e}
thus avoiding comparing the statistic pertaining to gene $k$ to statistics pertaining to other genes.
Our value of $B$, $B=10000$, is much larger than Storey and Tibshirani's $B=100$.

We can also introduce a simplified version of \eqref{eq:permutation_e}:
\begin{equation}\label{eq:permutation_e_simple}
  e_k
  :=
  \frac{T_k}
  {
    \frac{1}{B}
    \sum_{b=1}^B
    T_k^{0b}
  },
\end{equation}
in analogy with \eqref{eq:ST}.
This version may be more intuitive,
but it is only approximately valid for large $B$
and ceases to be valid for small $B$.
The difference shows, e.g., in the fact that $e_k$ defined via \eqref{eq:permutation_e} is bounded above by $B+1$
whereas $e_k$ defined via \eqref{eq:permutation_e_simple} is potentially unbounded;
such a difference can be significant if $B$ is small.
When $B=10000$,
there is not much difference between using \eqref{eq:permutation_e} and using \eqref{eq:permutation_e_simple}:
see the right panel of Figure~\ref{fig:BRCA_10000},
which uses \eqref{eq:permutation_e_simple}.

\begin{figure}
  \begin{center}
    \includegraphics[width=0.48\textwidth]{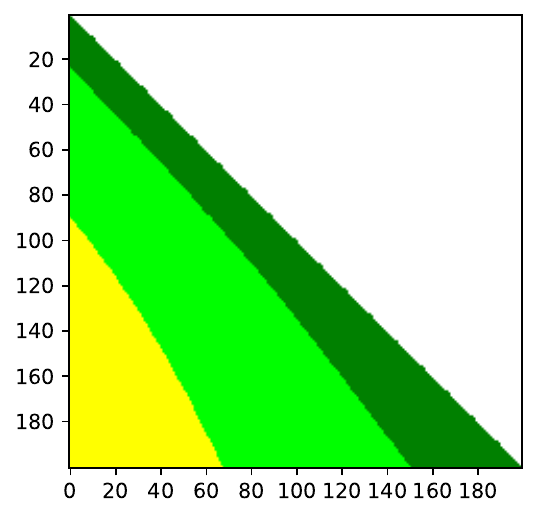}
    \includegraphics[width=0.48\textwidth]{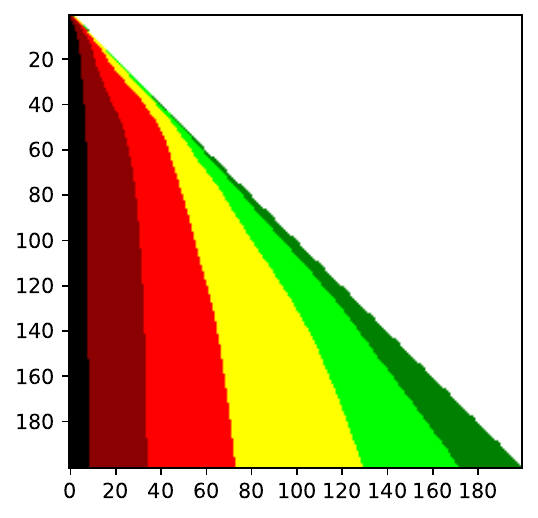}
  \end{center}
  \caption{Left panel:
    the top-left $200\times200$ corner of the arithmetic-mean discovery matrix
    for the \texttt{BRCA} dataset for $B:=100$.
    Right panel: its simplified version whose lack of validity is visible.}
  \label{fig:BRCA_100}
\end{figure}

A useful role of the version~\eqref{eq:permutation_e_simple} may be to check
whether the value of $B$ in \eqref{eq:permutation_e} is sufficiently large.
In the case of Figure~\ref{fig:BRCA_10000},
the approximation is good, which suggests that $B$ is sufficiently large.
However, in the case of Figure~\ref{fig:BRCA_100},
where $B=100$ (as in \citet{Storey/Tibshirani:2003}),
the right-hand panel, which uses \eqref{eq:permutation_e_simple},
looks far too good to be valid.
On the other hand, the left-hand panel,
which uses \eqref{eq:permutation_e},
is valid but extremely conservative.

\begin{figure}
  \begin{center}
    \includegraphics[width=0.48\textwidth]{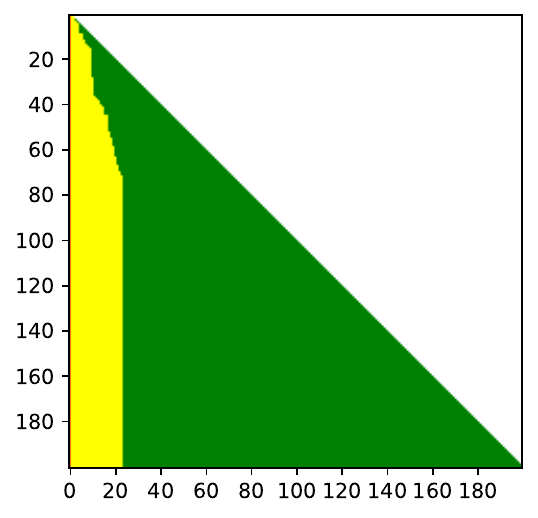}
    \includegraphics[width=0.48\textwidth]{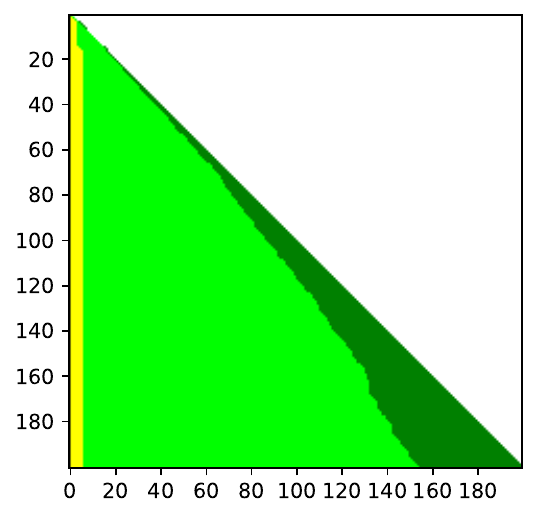}
  \end{center}
  \caption{Left panel:
    the top-left corner of the GWGS discovery \p-matrix for the \texttt{BRCA} dataset
    for Fisher's thresholds 1\% and 5\%, assuming independence.
    Right panel:
    analogous picture for Jeffreys's thresholds applied to the VS bounds for the entries of this matrix.}
  \label{fig:GWGS_emp_ind}
\end{figure}

Results given by \texttt{hommel} are either poor (when independence is assumed)
or extremely poor (under arbitrary dependence).
The former are given in Figure~\ref{fig:GWGS_emp_ind}
and the latter are given in Appendix~\ref{app:GWGS}
(Figure~\ref{fig:GWGS_emp_gen}).

\begin{table}
  \caption{The Benjamini--Hochberg and Benjamini--Yekutieli procedures
    applied to the \texttt{BRCA} dataset
    (the three entries of ``1'' are unreliable as they are based on a zero \p-value).
  \label{tab:emp_BH}}
  \begin{center}
    \begin{tabular}{lrr}
      assumption & $5\%$ & $1\%$ \\
      \hline
      independence & 88 & 1 \\
      arbitrary dependence & 1 & 1
    \end{tabular}
  \end{center}
\end{table}

The Benjamini--Hochberg procedure
\citep{Benjamini/Hochberg:1995}
rejects 88 null hypotheses at FDR $q:=0.05$
and 1 null hypothesis at FDR $q:=0.01$ for Storey and Tibshirani's list of \p-values
(of course, we will always reject at least 1 null hypothesis
because of the zero \p-value on their list).
However, this procedure assumes independence.
Under arbitrary dependence,
we can control FDR by replacing $q$ by $q/\sum_{k=1}^K k^{-1}$
\citep[Theorem 1.3]{Benjamini/Yekutieli:2001}.
This leads to rejecting 1 null hypothesis even at FDR $0.05$,
which is both poor and unwarranted.
These results are summarized in Table~\ref{tab:emp_BH}.

\section{Efficient implementation of Algorithm~\ref{alg:D} for the arithmetic mean}
\label{app:AM}

Algorithm~\ref{alg:D} is a generic algorithm that works for any symmetric \e-merging function $F$.
In general, computing one row of the discovery \e-matrix takes time $O(K^3)$
if we assume that the base \e-merging function $F$ can be computed in time linear in the number of arguments.
This assumption is correct for the arithmetic mean and, provided the arguments are sorted, the Simes \e-merging function
(see Appendix \ref{app:general}).
The overall computational complexity for the full discovery \e-matrix is very high, $O(K^4)$.

\begin{algorithm}[bt]
  \caption{Arithmetic-mean discovery matrix $\AM$}
  \label{alg:AM}
  \begin{algorithmic}[1]
    \Require
      An increasing sequence of \e-values $e_1\le\dots\le e_K$.
    \State $s_1 := e_1$\label{ln:preprocessing_1_start}
    \For{$k=2,\dots,K$}
      \State $s_k := s_{k-1}+e_k$\label{ln:preprocessing_1_end}
    \EndFor
    \For{$r=1,\dots,K$}\label{ln:preprocessing_2_start}
      \State $\sigma_{r,r-1} := e_{K-r+1}$
      \For{$j=r-2,\dots,0$}\label{ln:tricky_loop}
        \State $\sigma_{r,j} := \sigma_{r,j+1} + e_{K-j}$\label{ln:preprocessing_2_end}
      \EndFor
    \EndFor
    \For{$r=1,\dots,K$}
      \For{$j=0,\dots,r-1$}
        \State $\AM_{r,j} := \sigma_{r,j} / (r-j)$
        \For{$i=1,\dots,K-r$}
          \State $e := (\sigma_{r,j}+s_i) / (r-j+i)$
          \If{$e < \AM_{r,j}$}
            \State $\AM_{r,j} := e$
          \EndIf
        \EndFor
      \EndFor
    \EndFor
  \end{algorithmic}
\end{algorithm}

A more efficient implementation of Algorithm~\ref{alg:D} for the arithmetic mean
is given as Algorithm~\ref{alg:AM},
which uses arrays $s_k$ (the sum of the first $k$ base \e-values)
and $\sigma_{r,j}$ (the sum of the base \e-values with indices in $S_{r,j}$ in the notation of Algorithm~\ref{alg:D}).
There is a preprocessing stage (lines~\ref{ln:preprocessing_1_start}--\ref{ln:preprocessing_1_end})
taking time $O(K)$
and another preprocessing stage (lines~\ref{ln:preprocessing_2_start}--\ref{ln:preprocessing_2_end})
taking time $O(K^2)$;
the loop in lines \ref{ln:tricky_loop}--\ref{ln:preprocessing_2_end}
is executed in the decreasing order of $j$,
and in particular it is not executed when $r=1$
(this also applies to two similar loops in Algorithm~\ref{alg:AM_plus}).
After that computing each row of the arithmetic mean discovery matrix takes time $O(K^2)$.
The overall time is $O(K^3)$.

\begin{algorithm}[bt]
  \caption{One row of the arithmetic-mean discovery matrix in time $O(K)$}
  \label{alg:AM_plus}
  \begin{algorithmic}[1]
    \Require
      Increasing sequence of \e-values $e_1\le\dots\le e_K$
      and row number $r\in\{1,\dots,K\}$ of the discovery matrix.
    \State $s_0 := 0$
    \For{$k=1,\dots,K-r$}
      \State $s_k := s_{k-1}+e_k$
    \EndFor
    \State $\sigma_{r-1} := e_{K-r+1}$
    \For{$j=r-2,\dots,0$}
      \State $\sigma_{j} := \sigma_{j+1} + e_{K-j}$
    \EndFor
    \State $k := K-r$\label{ln:k_init}
    \For{$j=0,\dots,r-1$}
      \State $\texttt{slope} := \frac{s_k+\sigma_j}{k+r-j}$\label{ln:it_start}
      \For{$i=k-1,\dots,0$}
        \State $\texttt{new\_slope} := \frac{s_i+\sigma_j}{i+r-j}$
        \If{$\texttt{new\_slope} > \texttt{slope}$} \textbf{break}\label{ln:break}
        \EndIf
	\State $k := i$\label{ln:k}
	\State $\texttt{slope} := \texttt{new\_slope}$
      \EndFor
      \State $\AM_{r,j} := \texttt{slope}$\label{ln:it_end}
    \EndFor
  \end{algorithmic}
\end{algorithm}

\begin{figure}
  \begin{center}
    \includegraphics[width=0.55\textwidth]{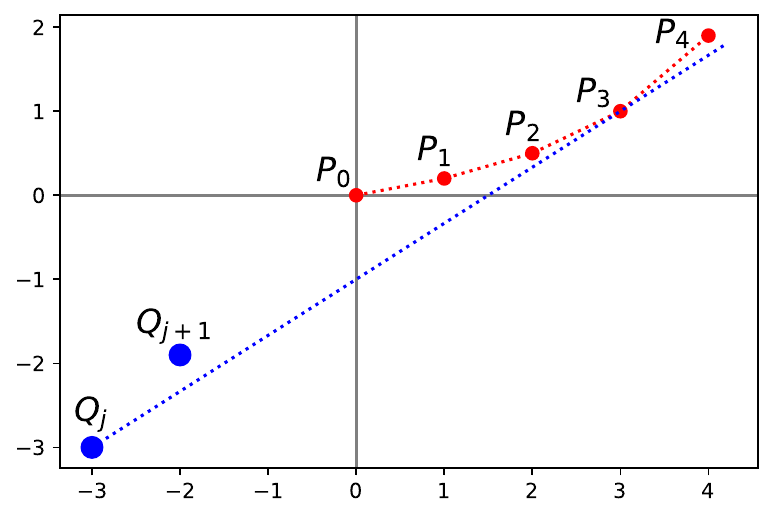}
  \end{center}
  \caption{Geometry behind Algorithm~\ref{alg:AM_plus}.}
  \label{fig:algorithm}
\end{figure}

An even more efficient implementation of Algorithm~\ref{alg:D} is given as Algorithm~\ref{alg:AM_plus}.
This algorithm computes one row of the arithmetic mean discovery matrix in time $O(K)$,
which gives the overall time $O(K^2)$.
Both $O(K)$ and $O(K^2)$ are clearly optimal in this context.
The ability to compute efficiently individual rows is useful when the discovery matrix is big;
e.g., it can be too big to fit in computer memory.

We are using essentially the same array $s$ as in Algorithm~\ref{alg:AM}
(now we extend it by adding $s_0:=0$),
and the array $\sigma$ in Algorithm~\ref{alg:AM_plus}
is one row of the array $\sigma$ in Algorithm~\ref{alg:AM};
\begin{align*}
  s_k &= e_1+\dots+e_k, & k&=0,\dots,K-r,\\
  \sigma_j &:= e_{K-r+1}+\dots+e_{K-j}, & j&=0,\dots,r-1.
\end{align*}
Of course, there is no need to recompute the array $s$ for each row $r$ of the discovery matrix.

The geometry behind Algorithm~\ref{alg:AM_plus} is shown in Figure~\ref{fig:algorithm}.
The coordinates of each of the points $P_k$, $k=0,\dots,K-r$, are $(k,s_k)$,
and the coordinates of the point $Q_j$, where $j\in\{0,\dots,r-1\}$, are
$-(r-j,\sigma_j)$.
Since the sequence $e_1,\dots,e_{K-r}$ is increasing,
connecting the points $P_0,P_1,\dots$ in this order
(see the red line in Figure~\ref{fig:algorithm})
gives us the graph of a convex function.

The command \textbf{break} in line~\ref{ln:break} means leaving the loop,
as in Python or \textsf{R};
in this context, it is equivalent to ``go to line~\ref{ln:it_end}''.
The variable $k$ in line \ref{ln:k} is the index of the ``current vertex'' $P_k$;
we start from the rightmost $P_k$ in line~\ref{ln:k_init} and then keep moving left.
Figure~\ref{fig:algorithm} illustrates the execution of Algorithm~\ref{alg:AM_plus}
when $k=3$, so that the current vertex is $P_3$.

For each $j=1,\dots,r$, the iteration of the loop in lines \ref{ln:it_start}--\ref{ln:it_end} of Algorithm~\ref{alg:AM_plus}
computes the slope of the straight line (shown in blue) passing through $Q_j$ and touching the red line from below.
The validity of the algorithm follows from the point $Q_{j+1}$ lying at or above the blue line,
for each $j=0,\dots,r-2$.
Let us check the last statement.
If $k<K-r$, the slope of the blue line is at most $e_{k+1}\le e_{K-r}$.
On the other hand, the slope of the line going from $Q_j$ to $Q_{j+1}$
is $e_{K-j}\ge e_{K-r}$.
It remains to consider the case $k=K-r$.
In this case, it suffices to notice that the slope $e_{K-j}$ of the line going from $Q_j$ to $Q_{j+1}$
is greater than or equal to the average of $e_1,\dots,e_{K-j}$,
which is the slope of the line going from $Q_j$ to $P_{K-r}$.

If we are only interested in the positions where discovery vectors or matrices exceed a given threshold,
we can also use algorithms described in \citet{Tian/etal:2021}.

\section{Conclusion}
\label{sec:conclusion}

The main technical tool of this paper, \e-values,
has important advantages over \p-values.
The advantage that we have found most useful here
is the easiness of merging \e-values:
the arithmetic average of \e-values is an \e-value,
and this is the only useful symmetric method of merging \e-values.
Other advantages were mentioned in Section~\ref{sec:introduction},
such as the open nature of \e-values allowing their sequential updating.

We have described methods for multiple hypothesis testing using \e-values
and demonstrated their use in simulation and empirical studies.
We believe that these methods, being simpler and more powerful,
are preferred to methods using \p-values
unless the final result must be stated in terms of \p-values.
Besides, our methods do not depend on the base \e-values being independent,
and under arbitrary dependence,
they are sometimes competitive with results based on \p-values
even when the final result is to be stated in terms of \p-values.

One of the obvious directions of further research is to extend our methods
to non-symmetric problems of multiple hypothesis testing (cf.\ \citet{Genovese/etal:2006}),
in which different \e-values may be assigned different weights.
Our procedure for multiple hypothesis testing is generic
and does not have to rely on unweighted arithmetic averaging.

\subsection*{Acknowledgments}

We are grateful to Peter Westfall for his advice about the literature
on Bayesian two-sample \textit{t}-tests.
We thank Glenn Shafer, Aaditya Ramdas,
and participants in the course ``Game-theoretic statistics'' (January--April 2021)
for helpful comments.
The presentation was greatly improved as result of the comments
by two referees, an Associate Editor, and the Editor (Sonia Petrone).
For most of our simulation and empirical studies in Sections~\ref{sec:simulation}--\ref{sec:empirical}
we used Python.
We also used the \textsf{R} package \texttt{hommel} \citep{Goeman/etal:2019R}
and a dataset available in the \textsf{R} package \texttt{qvalue} \citep{Storey/etal:2019}.

V.~Vovk's research has been partially supported by Amazon, Astra Zeneca, and Stena Line.
R.~Wang is supported by  the Natural Sciences and Engineering Research Council of Canada 
(RGPIN-2018-03823, RGPAS-2018-522590).

\appendix

\section{More \e-variables for testing $N(0,1)$}
\label{app:sigma}

\begin{figure}[b]
  \begin{center}
    \includegraphics[width=0.48\textwidth]{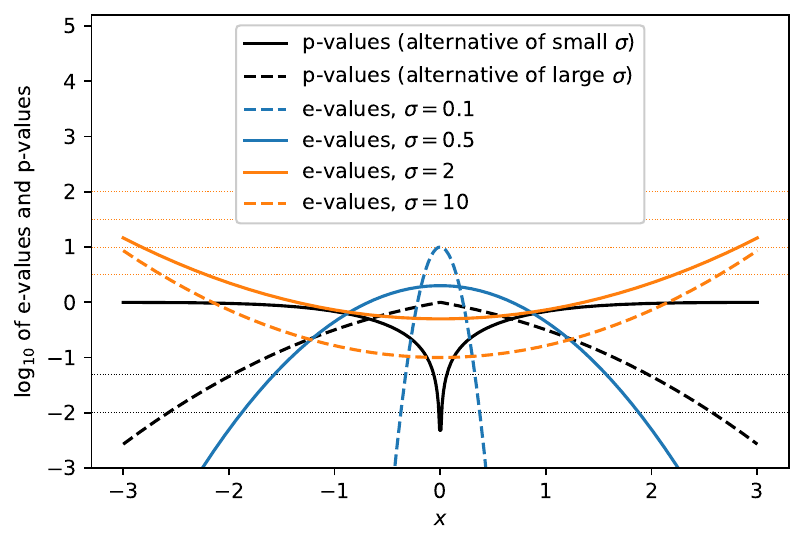}
  \end{center}
  \caption{Some \p-values and \e-values
    for testing $N(0,1)$ embedded into the family $N(0,\sigma^2)$.
    This complements the right panel of Figure~\ref{fig:toy}.}
  \label{fig:toy_sigma}
\end{figure}

Embedding the null hypothesis $N(0,1)$ into the family $N(\mu,1)$,
as in the right panel of Figure~\ref{fig:toy},
is not the only option,
and Figure \ref{fig:toy_sigma} gives results
for the family $N(0,\sigma^2)$, $\sigma>0$ being the standard deviation.
The figure shows the likelihood ratios $\dd N(0,\sigma^2)/\dd N(0,1)$
as \e-variables for a range of $\sigma$,
and it shows $P_<(x):=\chi^2(x^2)$ and $P_>(x):=1-\chi^2(x^2)$
as \p-variables,
where $\chi^2$ is the distribution function of $\chi^2$ with 1 degree of freedom
(with $P_<$ based on $\sigma<1$ as the alternative hypothesis,
and $P_>$ based on $\sigma>1$).

To get a non-trivial two-sided confidence interval for $\sigma$,
we need to merge the two \p-variables
(by, say, using the Bonferroni merging function
$B(P_<,P_>):=2\min(P_<,P_>)$)
or merge two of the \e-variables, one for $\sigma>1$ and the other for $\sigma<1$
(which can be done more efficiently, simply by averaging them).

\section{Using other \e-merging functions}
\label{app:general}

In this appendix, we briefly explore discovery \e-matrices using \e-merging functions $F$ other than the arithmetic mean.
For general \e-merging functions $F$
(not necessarily satisfying \eqref{eq:F-mono}),
we consider \emph{regularized} discovery \e-matrices
\[
  \bar D_{r,j}
  :=
  \min_{j'\le j}
  D_{r,j'}.
\]
They are lower bounds on $\eP^{g_R}(\{1,\dots,j\}\mid\omega)$.
Algorithm~\ref{alg:D} can be adapted to produce the regularized discovery matrix
by replacing $K-r$ with $K$ in line~\ref{ln:K-r}.

The \emph{Bonferroni \e-merging function} is the following lower bound for \eqref{eq:AM}:
\begin{equation}\label{eq:lower_1}
  B(e_1,\dots,e_n)
  :=
  \frac{1}{n}
  \max_{i\in\{1,\dots,n\}}
  e_i.
\end{equation}
A better lower bound for \eqref{eq:AM} is given by the \emph{Simes \e-merging function}
\begin{equation}\label{eq:lower_2}
  S(e_1,\dots,e_n)
  :=
  \max_{i\in\{1,\dots,n\}}
  \frac{i e_{[i]}}{n},
\end{equation}
where $e_{[i]}$ is the $i$th largest \e-value among $e_i$, $i\in\{1,\dots,n\}$
\citep[end of Section~6]{Vovk/Wang:2021}:
$e_{[1]},\dots,e_{[n]}$ is the permutation of $e_{1},\dots,e_{n}$ satisfying $e_{[1]}\ge\dots\ge e_{[n]}$.
While \eqref{eq:lower_1} sometimes violates \eqref{eq:F-mono}
(consider, e.g., the case $\mathbf{e}=(e,\dots,e)$ in \eqref{eq:F-mono}),
\eqref{eq:lower_2} is guaranteed to satisfy it.

\begin{proposition}
  For all $\mathbf{e}\in[0,\infty]^*$ and $e\ge\max(\mathbf{e})$,
  \begin{equation}\label{eq:S-mono}
    S(\mathbf{e},e) \ge S(\mathbf{e}).
  \end{equation}
\end{proposition}

\begin{proof}
  Without loss of generality, suppose $\mathbf{e}=(e_1,\dots,e_n)$
  with $e_n\le\dots\le e_1\le e$.
  Let $i$ be such that $S(\mathbf{e}) = i e_{i} / n$.
  By the definition of $S(\mathbf{e},e)$,
  \[
    S(\mathbf{e},e)
    \ge
    \frac{i+1}{n+1} e_{i}
    \ge
    \frac{i}{n} e_{i}
    =
    S(\mathbf{e}).
    \qedhere
  \]
\end{proof}

\noindent
Therefore, regularization is never needed for the discovery \e-matrices
based on the Simes \e-merging function.

Our discussion of \e-values and \p-values in Section~\ref{sec:comparison}
suggests that the function $t\mapsto 1/t$ transforms \e-values into \p-values
(cf.\ \eqref{eq:e-to-p_calibrator})
and transforms \p-values into approximate \e-values
(cf.\ \eqref{eq:calibrator} for a small $\kappa\in(0,1)$);
of course, the word ``approximate'' is used here in a crude sense
(in the spirit of the algorithmic theory of randomness).
Under this correspondence, the Bonferroni \e-merging function \eqref{eq:lower_1}
turns into the Bonferroni merging function for \p-values,
and the Simes \e-merging function \eqref{eq:lower_2}
turns into the Simes merging function for \p-values.
The dominating arithmetic-mean \e-merging function \eqref{eq:AM}
corresponds to using the harmonic mean for merging \p-values,
and indeed the harmonic mean has been discussed recently in this role \citep{Wilson:2019},
sometimes with a similar justification based on the VS bound \citep{Held:2019}.
However, the harmonic mean is not a valid function for merging \p-values \citep{Goeman/etal:2019}
unless multiplied by, say, $2.5 \ln K$ for $K\ge3$ \citep{Vovk/Wang:2020}.

\begin{algorithm}[b]
  \caption{Regularized Bonferroni discovery \e-matrix $\BM$}
  \label{alg:BM}
  \begin{algorithmic}[1]
    \Require
      An increasing sequence of \e-values $e_1\le\dots\le e_K$.
    \State{$a:=\infty$}\label{ln:a}
    \For{$j=0,\dots,K-1$}
      \State{$B := e_{K-j}/(K-j)$}\label{ln:B}
      \If{$a > B$}
        {$a := B$}\label{ln:adjust}
      \EndIf
      \For{$r=j+1,\dots,K$}
        \State{$\BM_{r,j} := a$}\label{ln:to-change}
      \EndFor
    \EndFor
  \end{algorithmic}
\end{algorithm}

With $F_{\mathbf{e}}(I)$ of \eqref{eq:F} specialized to the Bonferroni lower bound
\begin{equation*}
  F_{\mathbf{e}}(I)
  =
  B_{\mathbf{e}}(I)
  :=
  \frac{1}{\left|I\right|}
  \max_{i\in I}
  e_i,
  \quad
  I\subseteq\{1,\dots,K\},
  \enspace
  I\ne\emptyset,
\end{equation*}
Algorithms~\ref{alg:R} and~\ref{alg:D} have the same interpretation as before
(although the results are not as good since they are based on more conservative \e-values).
However, they simplify, especially Algorithm~\ref{alg:D},
whose regularized Bonferroni implementation is given as Algorithm~\ref{alg:BM}.
In line~\ref{ln:a} of Algorithm~\ref{alg:BM}
we initialize the adjusted Bonferroni \e-value,
in line~\ref{ln:B} we compute the raw Bonferroni \e-value,
and in line~\ref{ln:adjust} we adjust it.
The algorithm produces a matrix $\BM$ with constant columns
and takes time $O(K)$ per column;
this time is spent simply by writing one value repeatedly.
The resulting computational complexity $O(K^2)$ is clearly the optimal one.

\begin{figure}
  \begin{center}
    \includegraphics[width=0.48\textwidth]{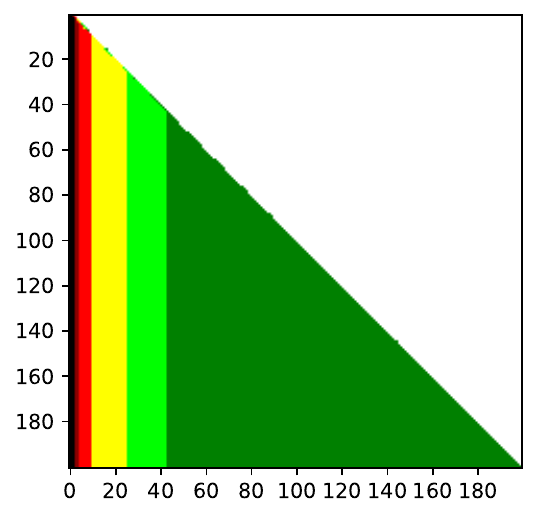}
    \includegraphics[width=0.48\textwidth]{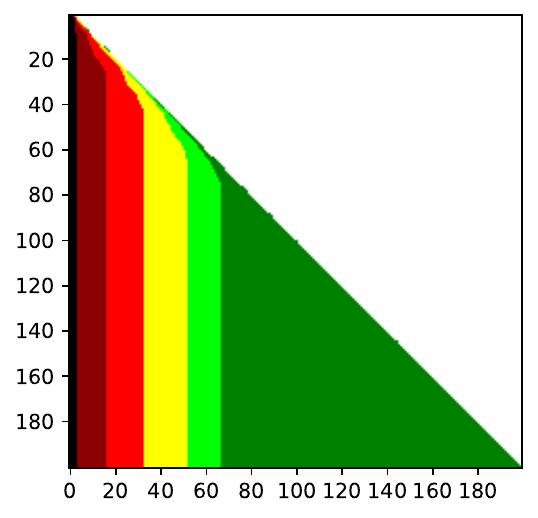}
  \end{center}
  \caption{Left panel: the regularized Bonferroni discovery \e-matrix (given by Algorithm~\ref{alg:BM})
    in the situation of Figure~\ref{fig:big},
    using Jeffreys's thresholds.
    Right panel: the corresponding Simes discovery \e-matrix.}
  \label{fig:BS}
\end{figure}

Whereas Bonferroni-type procedures often perform well when the goal is family-wise validity
(see, e.g., \citet[Figures~3 and~4]{Vovk/Wang:2021}),
their performance tends to deteriorate for less demanding notions of validity.
(In terms of \p-values, this phenomenon is discussed in, e.g., \citet[Section 1]{Goeman/Solari:2011rejoinder}.)
Comparing the left panel of Figure~\ref{fig:BS} with the upper left panel of Figure~\ref{fig:big}
we can see that the \e-Bonferroni method is much worse than arithmetic averaging
when the goal is to control the number of false discoveries.

The poor performance of the \e-Bonferroni method is clear already
from the upper left panel of Figure~\ref{fig:big}:
the areas of different colours are far from been vertical at the top, where they curve left.
It is clear that every discovery \e-matrix that is dominated by this one
and has vertical boundaries between different colours (such as \e-Bonferroni)
is going to be much worse.

The Bonferroni discovery \e-matrix (unregularized)
can be obtained by replacing $a$ with $B$ in line~\ref{ln:to-change} of Algorithm~\ref{alg:BM}
(and optionally ignoring all other lines containing $a$).
The left panel of Figure \ref{fig:BS} does not change if we remove the regularization.

The right panel of Figure~\ref{fig:BS} shows the Simes \e-matrix,
based on \eqref{eq:lower_2},
in the situation of Figure~\ref{fig:big}.
It is intermediate between $\AM$ and Bonferroni and, remarkably,
it looks better than the GWGS discovery \p-matrix
transformed by applying the VS bound to its elements
(the lower right panel of Figure~\ref{fig:big}).

\begin{remark}
  We can quantify the quality of the lower bounds \eqref{eq:lower_1} and \eqref{eq:lower_2}
  of the arithmetic mean $F$ by the inequalities
  \begin{align*}
    &B(e_1,\dots,e_n)
    \le
    S(e_1,\dots,e_n)
    \le
    F(e_1,\dots,e_n),\\
    &1
    \le
    \frac{F(e_1,\dots,e_n)}{S(e_1,\dots,e_n)}
    \le
    \sum_{k=1}^{n} \frac1k
    \le
    \ln n + 1,\\
    &1
    \le
    \frac{F(e_1,\dots,e_n)}{B(e_1,\dots,e_n)}
    \le
    n,
    \quad
    1
    \le
    \frac{S(e_1,\dots,e_n)}{B(e_1,\dots,e_n)}
    \le
    n,
  \end{align*}
  all of which are tight (achievable as equality for any $n$),
  apart from ${}\le\ln n + 1$
  (which is tight only for $n=1$).
  Since $n:=\left|I\right|\le K$ in \eqref{eq:D_e_first},
  this gives bounds for the ratios of the corresponding elements of the discovery \e-matrices
  built on top of $B$, $S$, and $F$.
\end{remark}

\section{Comparison with the GWGS procedure}
\label{app:GWGS}

First we discuss the GWGS multiple testing procedure
in the form described in \citet[Section~2]{Goeman/Solari:2011local} and in terms of our definitions.
Let $F:\cup_{n=1}^{\infty}[0,1]^n\to[0,1]$ be a \emph{\p-merging function},
i.e., a monotonic function transforming \p-variables into a \p-variable:
whenever $P_1,\dots,P_n$ are \p-variables for some $n\in\{1,2,\dots\}$,
$F(P_1,\dots,P_n)$ is a \p-variable.
Suppose that $F$ is symmetric.
With such an $F$ we can associate the following regularized analogue of \eqref{eq:D_e_first} in terms of \p-values:
\begin{equation}\label{eq:D_p}
  \bar D_{\p,F}^R(j)
  :=
  \max_{I:\left|R\setminus I\right|\le j}
  F(p_i, i\in I)
  \ge
  \pN^{g_R}
  (\{j+1,j+2,\dots\}\mid\omega),
\end{equation}
where the \p-test $P$ is defined by $P_\theta:=F(P_k:k\in I_{\theta})$
(analogously to \eqref{eq:e-test});
we leave the dependence on $p_1,\dots,p_K$ implicit,
following Goeman and Solari and similarly to the case of \e-values.

Goeman and Solari prefer a kind of inverse to the function \eqref{eq:D_p},
which they denote $f_{\alpha}(R)$,
suppressing the dependence on $p_1,\dots,p_K$;
we consider it as function of $\alpha\in[0,1]$,
which is interpreted as significance level.
We will see that this function satisfies
\begin{equation}\label{eq:connection}
  f_{\alpha}(R) > j
  \Longleftrightarrow
  \bar D_{\p,F}^R(j)\le\alpha
\end{equation}
(and this equivalence can serve as definition of $f$).
Therefore, it gives us the same lower \p-confidence bound on the number of true discoveries
at significance level~$\alpha$.

For the reader familiar with \citet{Goeman/Solari:2011local},
we will check that their definition indeed satisfies \eqref{eq:connection}.
They first define their bound
\begin{equation*}
  t_{\alpha}(R)
  :=
  \max\{\left|I\right|\mid I\subseteq R, I\notin\XXX\}
\end{equation*}
on the number of false discoveries,
where
\[
  \XXX
  :=
  \{I \mid \forall J\supseteq I: J\in\UUU\}
\]
are the subsets of $\{1,\dots,K\}$ rejected by the closed testing procedure,
and
\[
  \UUU
  :=
  \{I \mid F(p_i,i\in I)\le\alpha\}
\]
are the subsets of $\{1,\dots,K\}$ rejected by $F$;
in general, $I$ and $J$ will run over the subsets of $\{1,\dots,K\}$.
Then they define their bound on the number of true discoveries as
\begin{equation}\label{eq:f_t}
  f_{\alpha}(R)
  :=
  \left|R\right|
  -
  t_{\alpha}(R).
\end{equation}
(And they refer to the true discoveries as false hypotheses
and to the false discoveries as true hypotheses.)
The equivalence \eqref{eq:connection} can be checked as follows:
\begin{align*}
  &f_{\alpha}(R) > j
  \Longleftrightarrow
  t_{\alpha}(R) < \left|R\right|-j\\
  &\Longleftrightarrow
  \max\{\left|I\right|\mid I\subseteq R, I\notin\XXX\} < \left|R\right|-j\\
  &\Longleftrightarrow
  \left(
    \forall I\subseteq R:
    I\notin\XXX \Rightarrow \left|I\right| < \left|R\right|-j
  \right)\\
  &\Longleftrightarrow
  \left(
    \forall I\subseteq R:
    \left|I\right| \ge \left|R\right|-j \Rightarrow I\in\XXX 
  \right)\\
  &\Longleftrightarrow
  \left(
    \forall I\subseteq R:
    \left|I\right| \ge \left|R\right|-j
    \Rightarrow
    (\forall J\supseteq I: J\in\UUU)
  \right)\\
  &\Longleftrightarrow
  \left(
    \forall J:
    \left|J\cap R\right| \ge \left|R\right|-j
    \Rightarrow
    J\in\UUU
  \right)\\
  &\Longleftrightarrow
  \left(
    \forall J:
    \left|R\setminus J\right| \le j
    \Rightarrow
    J\in\UUU
  \right)\\
  &\Longleftrightarrow
  \max_{J:\left|R\setminus J\right| \le j}
  F(p_i,i\in J) \le \alpha
  \Longleftrightarrow
  \bar D_{\p,F}^R(j) \le \alpha.
\end{align*}

In Section~\ref{sec:simulation} we mentioned that \citet[Section~3]{Goeman/Solari:2011rejoinder}
introduced a version of the notion of a discovery \p-vector.
Namely they introduced the confidence distribution whose quantile function is
$\alpha\mapsto t_{\alpha}(R)$.
We can interpret \eqref{eq:connection} as $\bar D^R_{\p,F}$
being the distribution function whose quantile function is
$\alpha\mapsto f_{\alpha}(R)$.
Since $f$ and $t$ are so closely connected (see \eqref{eq:f_t}),
the discovery \p-vector is closely connected
to Goeman and Solari's confidence distribution function.

\subsection*{A property of completeness for the GWGS procedure}

\citet{Goeman/etal:2021} have shown that the GWGS procedure
is the only admissible one for controlling true discoveries.
Doesn't this mean that the \e-version of this procedure, Algorithm~\ref{alg:R},
is inadmissible?

Similarly to \eqref{eq:D_e_first},
the interpretation of the property of validity \eqref{eq:D_p}
in terms of a Fisher-type disjunction is:
the rejection set $R$ contains more than $j$ true discoveries
unless the outcome $\omega$ is $\bar D_{\p,F}^R(j)$-strange.
We can indeed obtain a property of validity of the same kind
(i.e., in terms of \p-values) for the procedure of Algorithm~\ref{alg:R}.
Our interpretation of \eqref{eq:D_e_first}
(under the monotonicity in $j$,
as in Proposition~\ref{prop:monotone-2}\eqref{it:m_row})
was that the rejection set $R$ contains more than $j$ true discoveries
unless $\omega$ has an \e-value of $\bar D_{\e,F}^R(j)$ or more.
Applying the canonical \e-to-\p\ calibrator \eqref{eq:e-to-p_calibrator},
we can see that $R$ contains more than $j$ true discoveries
unless $\omega$ has a \p-value of $1/\bar D_{\e,F}^R(j)$ or less.
We have the same property of validity,
but with $1/\bar D_{\e,F}^R(j)$ in place of $\bar D_{\p,F}^R(j)$.
By Goeman et al.'s result,
the procedure with $1/\bar D_{\e,F}^R(j)$ is either a GWGS procedure
or inadmissible.
Since the operation of \e-to-\p\ calibration is so crude,
there is no doubt that this procedure is inadmissible in non-degenerate cases.
The source of its inadmissibility is the inefficiency of converting \e-values into \p-values,
and it can be shown that Algorithm~\ref{alg:R} itself is admissible
when $F$ is arithmetic averaging.

\subsection*{More results for the \texttt{hommel} package}

\begin{figure}
  \begin{center}
    \includegraphics[width=0.48\textwidth]{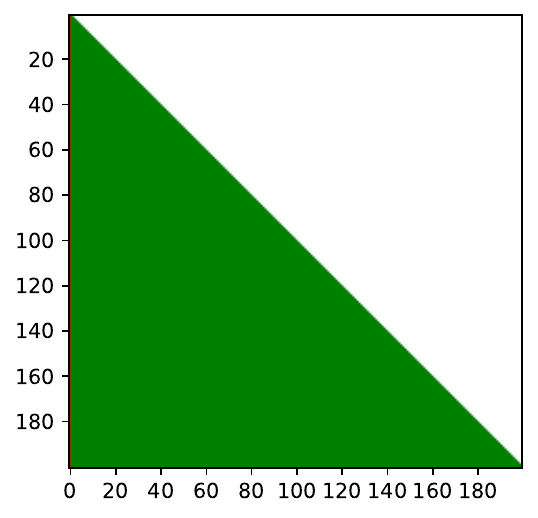}
    \includegraphics[width=0.48\textwidth]{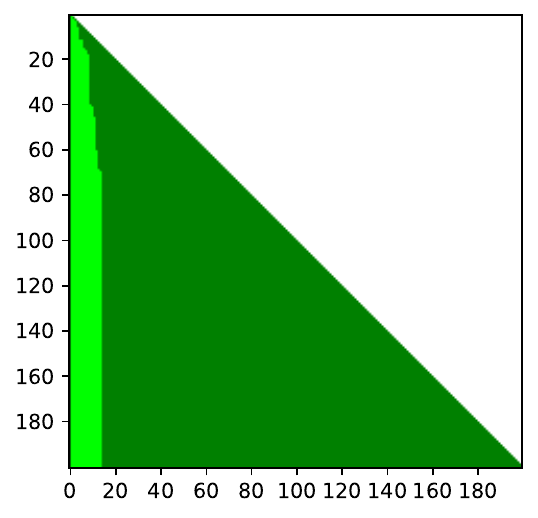}
  \end{center}
  \caption{Left panel:
    the top-left corner of the GWGS discovery \p-matrix for the \texttt{BRCA} dataset
    for Fisher's thresholds 1\% and 5\%, under arbitrary dependence.
    Right panel: analogous picture with each entry replaced by the corresponding VS bound
    and using Jeffreys's thresholds.}
  \label{fig:GWGS_emp_gen}
\end{figure}

In conclusion, we give one more figure demonstrating the work of the \texttt{hommel} package.
Figure~\ref{fig:GWGS_emp_gen} shows results (very poor) for the \texttt{BRCA} dataset without assuming independence.
In particular, the right panel is much worse than the left panel of Figure~\ref{fig:BRCA_10000},
which also does not assume independence.

\section{Generalized Bayes and boosting a weak signal}
\label{app:boosting}

When defining the base \e-values for use in our simulation studies
we just used the likelihood ratio $E(x)$ defined by \eqref{eq:E}.
This is the simplest version of a Bayes factor.
It usually works very well, but in some cases can be improved.
Later in this appendix we will see an example where a weak signal needs to be boosted,
but we start from developing tools that will allow us to do so.

Let us choose a constant $\eta>0$ (the \emph{learning rate})
and refer to
\begin{equation}\label{eq:E_gen}
  E_{\eta}(x)
  :=
  \frac{1}{c}
  E(x)^{\eta}
  =
  \frac{1}{c}
  \exp(\eta\delta x - \eta\delta^2/2)
\end{equation}
as the \emph{generalized Bayes factor}
(see, e.g.,
\citet[2.4]{Grunwald/Ommen:2017} and references therein).
Here $c>0$ is the normalizing constant ensuring $\int E_{\eta} \d N(0,1) = 1$;
a simple calculation gives
\[
  c
  =
  \exp
  \left(
    \eta (\eta-1) \delta^2/2
  \right).
\]
Plugging this into \eqref{eq:E_gen} we obtain
\begin{equation*}
  E_{\eta}(x)
  =
  \exp(\eta\delta x - \eta^2\delta^2/2).
\end{equation*}
This gives a useful interpretation of the generalized Bayes factor:
it is still the likelihood ratio,
but we replace the true alternative $N(\delta,1)$ by a false one, $N(\eta\delta,1)$.
For $\eta>1$ we are boosting the difference between the null and alternative hypotheses.

\begin{figure}
  \begin{center}
    \includegraphics[width=0.48\textwidth]{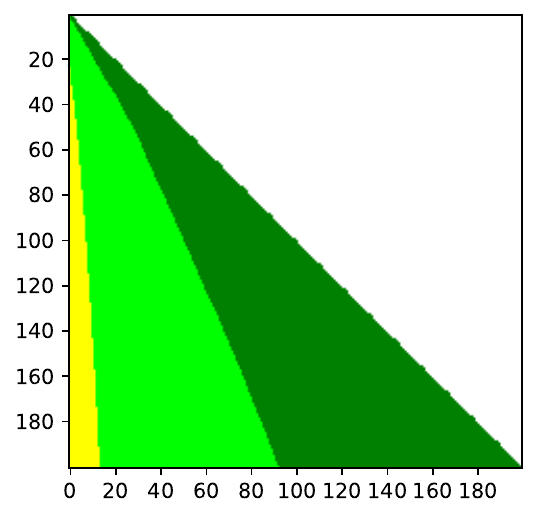}
    \includegraphics[width=0.48\textwidth]{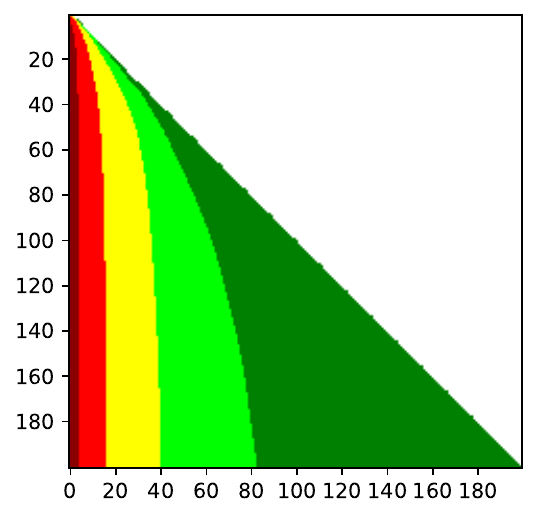}
  \end{center}
  \caption{Left panel:
    the top-left $200\times200$ corner of the arithmetic-mean discovery matrix
    for the simulation data with 10,000 observations, $10\%$ of false hypotheses, and weak signal,
    using Bayes factors as base \e-values,
    as described in text.
    Right panel: using generalized Bayes factors with learning rate $\eta=2$.
    Both panels use Jeffreys's thresholds.}
  \label{fig:sim_AM10000}
\end{figure}

\begin{figure}
  \begin{center}
    \includegraphics[width=0.48\textwidth]{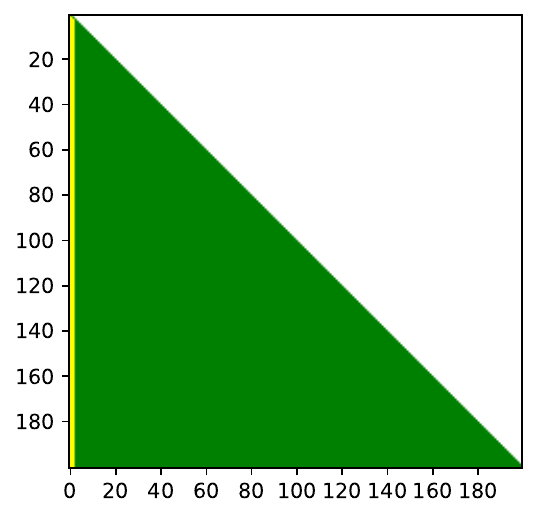}
    \includegraphics[width=0.48\textwidth]{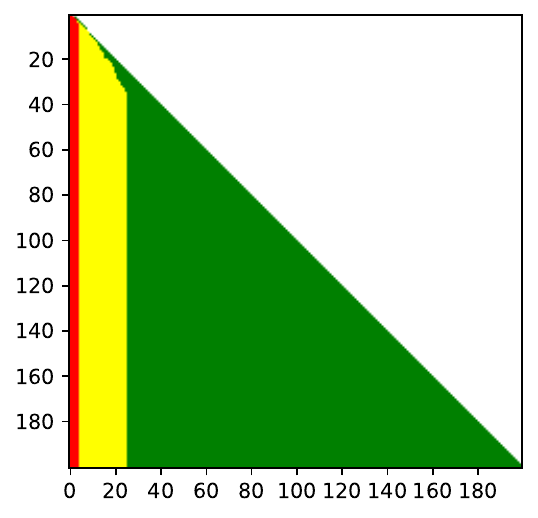}
  \end{center}
  \caption{Left panel:
    the top-left $200\times200$ corner of the GWGS discovery \p-matrix
    for the simulation data with 10,000 observations, $10\%$ of false hypotheses, and weak signal,
    under general dependence.
    Right panel: assuming independence.
    The colour code is based on Fisher's thresholds.}
  \label{fig:sim_GWGS10000}
\end{figure}

One situation in which the likelihood ratio \eqref{eq:E} does not work well
is where we have a large number of false null hypotheses,
but the true data-generating distributions are fairly close to the null hypotheses
(as it were, we have a weak signal).
Figures~\ref{fig:sim_AM10000}--\ref{fig:sim_GWGS10000} illustrate the case of 10,000 null hypotheses $N(0,1)$
of which 1000 are false, the true alternatives being $N(-2,1)$
(which makes the signal much weaker than in Section~\ref{sec:simulation}).
In the left panel of Figure~\ref{fig:sim_AM10000} we use the Bayes factor \eqref{eq:E},
whereas in its right panel we use the generalized Bayes factor \eqref{eq:E_gen} for $\eta=2$.
Using the generalized Bayes factor greatly improves the discovery \e-matrix.
The results for the GWGS procedure are given in Figure~\ref{fig:sim_GWGS10000};
they look poor, particularly so for arbitrary dependence.

\begin{table}
  \caption{The Benjamini--Hochberg and Benjamini--Yekutieli procedures
    applied to the simulation data with 10,000 observations, $10\%$ of false hypotheses, and weak signal
    for FDR $5\%$ and $1\%$.\label{tab:sim_BH_app}}
  \begin{center}
    \begin{tabular}{lrr}
      assumption & $5\%$ & $1\%$ \\
      \hline
      independence & 84 & 18 \\
      arbitrary dependence & 10 & 0
    \end{tabular}
  \end{center}
\end{table}

Table~\ref{tab:sim_BH_app} gives the numbers of rejections
for the Benjamini--Hochberg and Benjamini--Yekutieli procedures.
In view of Figure~\ref{fig:sim_AM10000} (right panel),
the results for arbitrary dependence are poor.

\begin{remark}
  The likelihood ratio of the true alternative to the null has well-known optimality properties
  as an \e-variable: see, e.g., \citet[\S2.2.1]{Shafer:2021}.
  Figure~\ref{fig:sim_AM10000} suggests that in multiple hypothesis testing
  the likelihood ratio may be far from being optimal.
  It would be interesting to explore this phenomenon theoretically.
\end{remark}

\section{Empirical study: ground truth}
\label{app:truth}

In Section~\ref{sec:empirical} we discussed a pioneering biomedical study
whose results were published a long time ago \citep{Hedenfalk/etal:2001}.
To evaluate the performance of various statistical techniques and their assumptions,
it is natural to analyze the developments in this area of biomedicine since 2001.

The main goal of \citet{Hedenfalk/etal:2001} was to test the hypothesis that different genes are expressed
by hereditary malignant breast tumours that are due to mutations in the BRCA1 and BRCA2 genes
and to identify differentially expressed genes.
A recent review \citep[Sections 2 and 5]{Wiggins/etal:2020} compares results of nine studies,
starting from \citet{Hedenfalk/etal:2001},
pursuing this goal and mostly using different biological samples
(therefore, not including \citet{Storey/Tibshirani:2003}).
The overlap between the lists of differentially expressed genes produced by different studies is poor.
In particular, only one gene has been identified as associated with BRCA1 by more than two studies.
This gene, TOB1, was among the genes identified in \citet[Figure 2A]{Hedenfalk/etal:2001}.
In Storey and Tibshirani's list of \p-values used in Figure~\ref{fig:GWGS_emp_ind} the TOB1 gene has rank 77;
in our list of \e-values used in Figure~\ref{fig:BRCA_10000} TOB1 has a slightly better rank of 51.

One reason \citep[Section 2]{Wiggins/etal:2020} for the poor overlap
between different studies is the genuine difficulty of the problem
of differentiating mutations in the two BRCA genes while controlling for potential confounders,
first of all the estrogen- and progesterone-receptor status and the subtype,
which are known to affect gene expression greatly.
For the dataset used in this paper,
differentiation between mutations in the two genes is facilitated, e.g.,
by all BRCA1 samples being negative for both estrogen and progesterone receptors
and majority of the BRCA2 samples being positive for both \citep[Table 1]{Hedenfalk/etal:2001}.
Other studies reported in \citet{Wiggins/etal:2020} tried to control for these confounders.

We can draw only limited conclusions from these follow-up studies.
There is often a big difference between the statistical null hypothesis
and the scientific hypothesis of interest.
Whereas there are genuine significant differences between the BRCA1 and BRCA2 samples in the dataset,
the differences are not necessarily due to their different BRCA status.
A possible lesson is that in our assumptions we should err on the side of caution
avoiding assuming independence or weak independence,
which lead us to expect very large numbers of discoveries.
\end{document}